\documentclass[a4paper]{article}
\usepackage[english]{babel}
\usepackage[utf8]{inputenc}
\usepackage{amsmath,amsthm}
\usepackage{amsfonts,amssymb}

\usepackage{tikz,tikz-cd} 

\usepackage{xifthen} 

\newcommand{\toposdefaut}{0}
\newcommand{\topos}[1][\toposdefaut]{ 
\ifthenelse{\equal{#1}{0}}{ \mathcal{T} }
{
\ifthenelse{\equal{#1}{1}}{ \mathcal{E} }{ #1 }
}
}

\newcommand{\set}{\text{Set}}

\newcommand{\spec}{\text{Spec }}

\newcommand{\scal}[2]{ \left\langle #1 , #2 \right\rangle }

\newcommand{\Q}{\mathbb{Q}}
\newcommand{\R}{\mathbb{R}}
\newcommand{\C}{\mathbb{C}}
\newcommand{\N}{\mathbb{N}}

\newcommand{\Ob}{\mathbb{O}}
\newcommand{\Char}{\mathbb{I}}

\newcommand{\Ecal}{\mathcal{E}} 
 
\newcommand{\Tcal}{\mathcal{T}}

\newcommand{\Gcal}{\mathcal{G}} 
\newcommand{\Hcal}{\mathcal{H}}

\newcommand{\Lcal}{\mathcal{L}}

\newcommand{\Ccal}{\mathcal{C}}

\usepackage{titlesec}
\titleformat{\subsubsection}[runin]{\normalfont}{\thesubsubsection}{0pt}{}[.]

\renewcommand{\thesubsubsection}{\arabic{section}.\arabic{subsubsection}}
\csname @addtoreset\endcsname{subsubsection}{section}

\newcommand{\block}[1]
{

\par \subsubsection{} #1

\bigskip}

\newcommand{\blockn}[1]{\par #1 \bigskip}

\newcommand{\Th}[1]
	{
	\bigskip	
	\textbf{Theorem : }{\itshape #1}
		
	\bigskip
	}

\newcommand{\Prop}[1]
	{

	\bigskip
	
	\textbf{Proposition : }{\itshape #1}
		
	\bigskip
	
	}

\newcommand{\Lem}[1]
	{

	\bigskip
	
	\textbf{Lemma : }{\itshape #1}
		
	\bigskip
	
	}

\newcommand{\Def}[1]
	{
	
	\bigskip
	
	\textbf{Definition : }{\itshape #1}
	
	\bigskip
	
	}

\newcommand{\Dem}[1]{
	
	\smallskip
	
	\textbf{Proof : } \par
	 {#1} $\square$
	 
	 \bigskip
}

\begin{document}

\pagestyle{plain}
\title{Toward a non-commutative Gelfand duality: Boolean locally separated toposes and Monoidal monotone complete $C^{*}$-categories}
\author{Simon Henry}

\maketitle

\begin{abstract}
** Draft Version ** To any boolean topos one can associate its category of internal Hilbert spaces, and if the topos is locally separated one can consider a full subcategory of square integrable Hilbert spaces. In both case it is a symmetric monoidal monotone complete $C^{*}$-category. We will prove that any boolean locally separated topos can be reconstructed as the classifying topos of ``non-degenerate" monoidal normal $*$-representations of both its category of internal Hilbert spaces and its category of square integrable Hilbert spaces. This suggest a possible extension of the usual Gelfand duality between a class of toposes (or more generally localic stacks or localic groupoids) and a class of symmetric monoidal $C^{*}$-categories yet to be discovered.
\end{abstract}

\renewcommand{\thefootnote}{\fnsymbol{footnote}} 
\footnotetext{\emph{email:} simon.henry@imj-prg.fr}
\renewcommand{\thefootnote}{\arabic{footnote}} 


\tableofcontents

\section{Introduction}

\blockn{\emph{This is a draft version. It will be replaced by a more definitive version within a few months. In the meantime any comments is welcome.}}

\blockn{This paper is part of a program (starting with the author's phd thesis) devoted to the study of the relation between topos theory and non-commutative geometry as two generalizations of topology. The central theme of this research project is the construction explained in section \ref{secPrelim} which naturally associate to any topos $\Tcal$ the $C^{*}$-category $\Hcal(\Tcal)$ of Hilbert bundles over $\Tcal$.}

\blockn{In the previous paper \cite{henry2014measure} we studied ``measure theory" of toposes and compare it to the theory of $W^{*}$-algebras (von Neumann algebras) through this construction of $\Hcal(\Tcal)$. At the end of the introduction of this previous paper one can find a table informally summing up some sort of partially dictionary between topos theory and operator algebra theory. The goal of the present work is somehow to provide a framework for making this dictionary a concrete mathematical result, by showing that a boolean locally separated topos $\Tcal$ can actually be completely reconstructed from any of the two symmetric monoidal\footnote{The term monoidal is taken with a slightly extended meaning: $\Hcal^{red}(\Tcal)$ has in general no unit object} $C^{*}$-categories $\Hcal(\Tcal)$ and $\Hcal^{red}(\Tcal)$ of Hilbert bundles and square-integrable Hilbert bundles.
}

\blockn{More precisely, we will show that if $\Tcal$ is a boolean locally separated topos, then $\Tcal$ is the classifying topos for non-degenerate\footnote{see definitions \ref{DefNondegenUnred} and \ref{DefRepReduced}.}  normal symmetric monoidal representations of either $\Hcal^{red}(\Tcal)$ and $\Hcal(\Tcal)$. These are not geometric theory, not even first order theory in fact, but we will prove an equivalence of categories between points of $\Tcal$ and these representations over an arbitrary base topos, see theorems \ref{mainThunred} and \ref{mainThRed} for the precise statement.
}

\blockn{Sections \ref{secPrelim},\ref{secMCcat_of_a_bool} and \ref{secLSTandSIHS} contain some preliminaries which are mostly, but not entirely, recall of previous work.}

\blockn{Section \ref{secStatement} contains the statement of the two main theorems of the present paper: theorem \ref{mainThunred}, which is the reconstruction theorem from the ``unreduced" category $\Hcal(\Tcal)$ and theorem \ref{mainThRed} which is the reconstruction theorem from the ``reduced" category $\Hcal^{red}(\Tcal)$ of square integrable Hilbert spaces.
}

\blockn{Sections \ref{secProofRed} and \ref{secProofUnred} contain respectively the proof of the ``reduced" theorem \ref{mainThRed} and the proof that the reduced theorem imply the unreduced theorem \ref{mainThunred}.

The key arguments, which actually use the specificity of the hypothesis ``boolean locally separated" seems to all be in subsection \ref{GeomorphonObjec}, and the rest of the proof seems to us, at least in comparison, more elementary (it is mostly about some computations of multilinear algebra) and also more general (the specific hypothesis are essentially\footnote{Lemma \ref{Fsurjection2} seems to be the unique exception} not used anymore). In fact, we have unsuccessfully tried to prove a result in this spirit for several years, and results of subsection \ref{GeomorphonObjec} have always been the main stumbling block.
}

\blockn{Finally we conclude this paper with section \ref{secToward}, where we explain how the results of this paper might maybe be extended into a duality between certain monoidal symmetric $C^{*}$-categories and certain geometric objects (presumably, localic groupoids or localic stacks) which would be a common extension of the usual Gelfand duality, the Gelfand duality for $W^{*}$-algebras, the Doplicher-Roberts reconstruction theorem for compact groups, and of course the results of the present paper. This would somehow constitute a sort of non-commutative Gelfand duality. Of course the existence of such a duality, and even its precise statement are at the present time highly conjectural, but we will try to highlight what are the main difficulties on the road toward such a result.

This conjectured duality is formally extremely similar to the reconstruction theorems obtained for algebraic stacks, as for example \cite{lurie2004tannaka}, \cite{brandenburg2014tensor}.

}

\section{General preliminaries}
\label{secPrelim}

\blockn{We will make an intensive use of the internal logic of toposes (i.e. the Kripke-Joyal semantics for intuitionist logic in toposes) in this paper. A reader unfamiliar with this technique can read for example sections $14$,$15$ and $16$ of \cite{mclarty1996elementary} which give a relatively short and clear account of the subject. Other possible references are \cite[chapter 6]{borceux3}, \cite[chapter VI]{maclane1992sheaves}, or \cite[D1 and D4]{sketches}.

}

\blockn{Because this paper is mostly about boolean toposes and monotone complete $C^{*}$-algebras we will assume that the base topos (that is the category of set) satisfies the law of excluded middle, but we won't need to assume the axiom of choice.

It is reasonable to think that these results can also be formulated and proved over a non-boolean basis, but the gain in doing so would be very small: over a non-boolean basis, any boolean topos and any monotone complete $C^{*}$-algebra is automatically defined over a boolean sub-locale of the terminal object, and hence can be dealt with in a boolean framework.

This being said, a large part of the proofs will take place internally in a non boolean topos and hence will have to avoid the law of excluded middle anyway.
 }

\blockn{We will also often have to juggle between internal and external logic or between the internal logic of two different toposes. We will generally precise what statement has to be interpreted ``internally in $\Tcal$" or ``externally", but we also would like to emphasis the fact that most of the time the context makes this completely clear: if an argument start by ``let $x \in X$" and that $X$ is not a set but an object of a topos $\Tcal$ it obviously means that we are working internally in $\Tcal$. This convention, is in fact completely similar with the usual use of context\footnote{``Context" is taken here in its formal ``type theoretic" meaning, i.e. the ``set" of all variable that has been declared at a given point of a proof.} in mathematics: when a mathematicians says something like ``let $x \in S$" (for $S$ a set) then what follows is actually mathematics internal to the topos $Sets_{/S}$ of sets over $S$, indeed everything being said after implicitly depends on a parameter $s \in S$, and if the conclusion actually does not depends on $s$ then it will be valid independently of the context only if $S$ was non-empty.

Our convention is hence that when we say something like ``let $x \in X$" where $X$ is an object of a topos $\Tcal$ then what follows is internal to the topos $\Tcal_{/X}$, or equivalently, internal to $\Tcal$ with a declared variable $x$, and this being true until the moment where this variable $x$ is ``removed" from the context (in classical mathematics, this is usually left implicit because, as soon as $X$ is non-empty, it is irrelevant, but in our case we will generally make it precise by saying that we are now working externally).

We will of course say explicitly in which topos we are working as soon as we think that it actually improve the readability, but we also think that this perspective makes (once we are used to it) the change from working internally into one topos from another topos as simple as introducing and forgetting abstract variables in usual mathematics and makes the text easier to read.

}

\blockn{All the toposes considered are Grothendieck toposes, in particular they all have a natural numbers object (see \cite[A2.5 and D5.1]{sketches}) ``$\N$" or ``$\N_{\Tcal}$" which is just the locally constant sheaf equal to $\N$.}

\blockn{A set, or an object $X$ of a topos, is said to be \emph{inhabited} if (internally) it satisfies $\exists x \in X$. For an object of a topos it corresponds to the fact that the map from $X$ to the terminal object is an epimorphism. }

\blockn{An object $X$ of a topos $\Tcal$ is said to be a bound of $\Tcal$ if subobjects of $X$ form a generating family of $\Tcal$ (i.e. is any object of $\Tcal$ admit a covering by subobject of $X$). Every Grothendieck topos admit a bound (for exemaple take the direct sum of all the representable sheaves for a given site of definition), in fact the existence of a bound together with the existence of small co-products characterize Grothendieck toposes among elementary toposes.}

\blockn{Let $\Tcal$ be an arbitrary topos. $\C_{\Tcal}$ is the object of ``continuous\footnote{also called Dedekind complex numbers} complex numbers" that is $\R_{\Tcal} \times \R_{\Tcal}$ where $\R_{\Tcal}$ is the object of continuous/Dedekind real numbers as defined for example in \cite[D4.7]{sketches}. In any topos (with a natural number object), $\C_{\Tcal}$ is a locale ring object.}

\blockn{When $P$ is a decidable proposition, i.e. if one does have $P$ or not-$P$, one denotes $\Char_P$ the real number defined by $\Char_P = 1$ if $P$ and $\Char_P$=0 otherwise. An object $X$ is said to be \emph{decidable} if for all $x$ and $y$ in $X$ the proposition $x = y$ is decidable, in which case we denote $\delta_{x,y}$ for $\Char_{x=y}$. i.e. $\delta_{x,y}$ is one if $x=y$ and zero otherwise. }

\blockn{By a ``Hilbert space of $\Tcal$", or a $\Tcal$-Hilbert space we mean an object $H$ of $\Tcal$, endowed with a $\C_{\Tcal}$-module structure and a scalar product $H \times H \rightarrow \C_{\Tcal}$ (linear in the second variable and anti-linear in the first), which satisfies internally all the usual axioms for being a Hilbert space, completeness being interpreted in term of Cauchy filters, or equivalently Cauchy approximations but not Cauchy sequences.
}

\blockn{$\Hcal(\Tcal)$ denotes the $C^{*}$-category\footnote{see for exemple \cite{wstarcat} for the definition of $C^{*}$-category.} whose objects are Hilbert spaces of $\Tcal$ and whose morphisms are ``globally bounded operators", that is linear maps $f : H \rightarrow H'$ which admit an adjoint and such that it exists an external number $K$ satisfying (internally) for all $x \in H$, $\Vert f(x) \Vert \leqslant K \Vert x \Vert $. The norm $\Vert f \Vert_{\infty}$ is then the smallest such constant $K$ (if we were not assuming the law of excluded middle in the base topos it would be an upper semi-continuous real number), the addition and composition of operators is defined internally , $f^{*}$ is the adjoint of $f$ internally, and this form a $C^{*}$-category.
}

\blockn{Because the tensor product of two Hilbert spaces can easily be defined, even in intuitionist mathematics, the category $\Hcal(\Tcal)$ is endowed with a symmetric monoidal structure.}

\blockn{Precise definition of a symmetric monoidal category, a symmetric monoidal functor, and a symmetric monoidal natural transformation can be found in S.MacLane's category theory books\footnote{Only in the second edition (1998).} \cite{maclane1998categories}, chapter \textsc{XI}, section $1$ and $2$.

Briefly, a symmetric monoidal category is a category endowed with a bi-functor $\_ \otimes \_ $, a specific ``unit" object $e$ and isomorphisms $e \otimes A \simeq A \otimes e \simeq A$,  $(A \otimes B) \simeq (B \otimes A)$ and $(A \otimes B) \otimes C \simeq A \otimes ( B \otimes C)$ which are natural and satisfies certain coherence conditions. A symmetric monoidal functor (``braided strong monoidal functor" in MacLane's terminology) is a functor between symmetric monoidal categories with a natural isomorphism $F(A \otimes B) \simeq F(A) \otimes F(B)$ which has to satisfy a certain number of coherence and compatibility relations. Finally a symmetric monoidal transformation (the adjective symmetric is actually irrelevant for natural transformation) is a natural transformation between symmetric monoidal functor which satisfy coherence conditions, stating that, up to the previously defined natural isomorphisms, $\eta_A \otimes \eta_b$ is the same as $\eta_{A \otimes B}$ and $\eta_{e}$ is the identity.
}

\blockn{Moreover, when we are talking about (symmetric) monoidal $C^{*}$-categories or symmetric monoidal $*$-functor between such categories we are always assuming that all the structural isomorphisms are in fact isometric isomorphisms, i.e. their inverse is their adjoint. For example, it is clearly the case for $\Hcal(\Tcal)$.}

\blockn{Finally, if $f:\Ecal \rightarrow \Tcal$ is a geometric morphism between two toposes, and $H$ is a Hilbert space of $\Tcal$ then $f^{*}(H)$ is a ``pre-Hilbert" space of $\Ecal$ but fails in general to be complete and separated, we denote by $ f^{\sharp}(H)$ its separated completion. $f^{\sharp}$ is a symmetric monoidal $*$-functor from $\Hcal(\Tcal)$ to $\Hcal(\Ecal)$. }

\section{Monotone complete $C^{*}$-categories and boolean toposes}
\label{secMCcat_of_a_bool}

\blockn{We recall that a $C^{*}$-algebra is said to be \emph{monotone complete} if every bounded directed net of positive operators has a supremum. A positive linear map between two monotone complete $C^{*}$-algebras is said to be \emph{normal} if it preserves supremum of bounded directed set of positive operators.

The theory of monotone complete $C^{*}$-algebras is extremely close to the theory of $W^{*}$-algebras, in fact it is well know that a monotone complete $C^{*}$-algebra having enough normal positive linear form is a $W^{*}$-algebra (see \cite[Theorem 3.16]{takesaki2003theoryI}). }

\blockn{When $\Tcal$ is a boolean topos, $\Hcal(\topos)$ is a monotone complete $C^{*}$-category in the sense that it has bi-products and the $C^{*}$-algebra of endomorphisms of any object is monotone complete. Indeed, because $\topos$ is boolean, the supremum of a bounded net of operators can be computed internally, and as the supremum is unique it ``patches up" into an externally defined map, see \cite[section 2]{henry2014measure}.

Monotone complete $C^{*}$-categories are extremely close to the $W^{*}$-categories studied in \cite{wstarcat}, in fact most of the result of \cite{wstarcat} which does not involve the existence of normal states (or the modular time evolution) also hold for monotone complete $C^{*}$-categories. We will review some of these results:
}

\blockn{If $C$ is a monotone complete $C^{*}$-category then we define the center $Z(C)$ of $C$ as being the commutative monotone complete $C^{*}$-algebra of endomorphisms of the identity functor of $C$. In the more general situation $Z(C)$ might fail to be a set and be a proper class, but we will not be concern by this issue because we proved in \cite[3.6]{henry2014measure} that $\Hcal(\topos)$ has a generator and hence, by results of \cite{wstarcat}, the algebra $Z(C)$ can be identified with the center of the algebra of endomorphisms of this generator. }

\blockn{If $A$ is an object of a monotone complete $C^{*}$-category then we define its central support $c(A) \in Z(C)$ by:

\[ c(A)_B := \sup_{f:A^{n} \rightarrow B \atop \Vert f \Vert < 1} f f^{*} \in Hom(B,B). \]

If the monotone complete category we are working with admit bi-product then $A^{n}$ denotes the bi-product of $n$-copies of $A$ and if not, then one can still make sense of a map $(f_1,\dots , f_n)$ from $A^{n}$ to $B$ as the data of $n$ maps from $A$ to $B$, $ff^{*}$ is $\sum f_i f_i ^{*}$ a,d $\Vert f \Vert$ is defined as $\Vert f f^{*}\Vert^{1/2}$. One can check that the supremum involved in the definition of $C(A)_B$ is directed by showing that the set of such ``$ff^{*}$" is in order preserving bijection with the set of $ff^{*}$ where $f$ is an arbitrary maps from $A^{n}$ to $B$ without condition on the norm, the bijection being obtained by multiplying $f$ by a convenient function of $ff^{*}$.

Equivalently, $c(A)$ can be defined as the smallest projection $c$ in $Z(C)$ such that $c_A=Id_A$, but we will need the fact that it is a directed supremum.

}

\blockn{One says that an object $A$ is \emph{quasi-contained} in an object $B$ if $c(A) \leqslant c(B)$. An object $A$ is said to be a \emph{generator} of a monotone complete $C^{*}$-category if and only if $c(A)=1$ i.e. if every other object is quasi-contained in $A$. }

\blockn{One can for example check that if two normal functors agree on a generator and its endomorphisms then they are isomorphic: it is an easy consequence of results of \cite{wstarcat} for $W^{*}$-categories and the proof can extended to monotone complete $C^{*}$-categories easily. }

\blockn{We conclude this section by briefly mentioning what quasi-containement mean in the case of $\Hcal(\Tcal)$: }

\block{\label{QcontaininHcal}\Prop{Let $\Tcal$ be a boolean topos and $H,H' \in \Hcal(\Tcal)$ two Hilbert spaces of $\Tcal$, then $H$ is weakly contained in $H'$ if and only if there exists a set $F$ of bounded operators from $H'$ to $H$ such that internally in $\Tcal$ the functions in $F$ spam a dense subspace of $H$.}

If this is true, we will say that $H$ is covered by the maps in $F$.

\Dem{If $H$ is weakly contained in $H'$ then $c(H) \leqslant c(H')$ hence $c(H')_H=1$. Rewriting this using the definition of $c(H')$ one gets:

\[ Id_H = \sup_{f:(H')^{n} \rightarrow H \atop \Vert f \Vert < 1} f f^{*} \in Hom(H,H) \]

But as we mentioned earlier, supremums of directed nets in $\Hcal(\Tcal)$ are computed internally, this means that this supremum converge internally for the strong operator topology. In particular, for any $h\in H$ one has $ff^{*}(h)$ which is arbitrarily close to $h$ when $f$ run through the (external) set:

\[ F_0 =\{f | f: (H')^{n} \rightarrow H, \Vert f \Vert < 1 \} \]

Hence taking $F$ to be the set of ``component" of maps in $F_0$, the sum of the images of maps in $F$ spam all of $H$.

\smallskip

Conversely, assume that $H$ is spammed by a family $F$ of external maps $f:H' \rightarrow H$. For any $f:H' \rightarrow H$, and in particular, for any $f \in F$ one has $c(H')_H \circ f = f \circ c(H')_{H'} = f$. The projector $c(H')_H$ is hence (internally in $\Tcal$) equal to the identity on the image of all the maps $f \in F$ and hence on all of $H$, i.e. $c(H')_H=Id_H$ which proves that $c(H) \leqslant c(H')$.
}

}






\section{Locally separated toposes and square integrable Hilbert spaces}
\label{secLSTandSIHS}
\blockn{We recall (see \cite[chapter II]{moerdijk2000proper}) that a topos $\topos$ is said to be \emph{separated} if its diagonal map $\topos \rightarrow \topos \times \topos$ (which is localic by \cite[B3.3.8]{sketches}) is proper, i.e. if, when seen as a $(\topos \times \topos)$-locale though the diagonal maps, $\topos$ is compact.
}

\blockn{In \cite[theorem 5.2]{henry2014measure} we proved that a boolean topos is separated if and only if it is generated by internally finite objects. One will use a slightly modified form of this result :}

\block{\label{SepImpFinitness}\Th{Let $\Tcal$ be a boolean separated topos, then $\Tcal$ admit a generating family of objects $(X)$ such that for each $X$ there exists an interger $n$ such that internally in $\Tcal$, the cardinal of $X$ is smaller than $n$.}
One will say that such objects are of \emph{bounded cardinal}.

\Dem{This theorem is an immediate consequence of \cite[theorem 5.2]{henry2014measure}: $\Tcal$ is generated by a familly of internally finite object $X$, but for each object $X$ of this generating familly and for each natural number $n$ one can define $X_n := \{x \in X | |X| \leqslant n\}$ whose cardinal is internally bounded by $n$, and as $X$ is internally finite the $(X_n)_{n \in \N}$ form a covering family of $X$ and hence the $(X_n)_{X,n}$ form a generating family fulfilling the property announced in the theorem.} }

\blockn{An object $X$ of boolean topos is said to be \emph{separating} if the slice topos $\topos_{/X}$ is separated. A boolean topos is said to be \emph{locally separated} if it admit an inhabited separating object, or equivalently if any object can be covered by separating objects, see \cite{henry2014measure} section 5 for more details. }

\blockn{If $X$ is an object of a boolean topos $\topos$ then one can define the $\Tcal$-Hilbert space $l^{2}(X)$ of square sumable sequences indexed by $X$. It can also be done in a non boolean topos but it require $X$ to be a decidable\footnote{In order to define the scalar product of two generators or to define the sum of a sequence we need that for any $x,y \in X$ $x = y$ or $x \neq y$.} object. Internally in $\topos$, the space $l^{2}(X)$ has generators $e_x$ for $x \in X$ such that $\scal{e_x}{e_y}=\delta_{x,y}$.}

\block{\label{l2weaklycontinl2}\Prop{Let $X$ be a bound of a boolean topos $\topos$ and $Y$ be a separating object of $\topos$ then $l^{2}(Y)$ is quasi-contained in $l^{2}(X)$.}

\Dem{Let $X$ be a bound and $Y$ any separating object. As $X$ is a bound, $Y$ can be covered by maps $h:U \rightarrow Y$ with $U \subset X$. Using the fact that $\Tcal_{/Y}$ is separated and boolean, we know that it is generated by objects of bounded cardinal and the image of a map whose domain is finite is also finite, and of smaller cardinal, hence $U$ admit a covering by sub-objects with bounded cardinal in $\Tcal_{/Y}$, hence we can freely assume that $U$ is itself of bounded cardinal in $\Tcal_{/Y}$.

One can then define (for each such map $h : U \rightarrow Y$) a map $\phi_h : l^{2}(X) \rightarrow l^{2}(Y)$ by $\phi_h(e_x)=0$ if $x \notin U$ and $\phi_h(e_x) = e_{h(x)}$ if $x \in U$. The fact that the object $U$ is finite with bounded cardinal over $Y$ mean that there exists an (external) integer $n$ such that each fiber of $h$ has cardinal smaller than $n$, this is exactly what we need to know to construct the adjoint of $\phi_h$ and to prove that $\phi_h$ is bounded (and hence extend into an operator).

Now as such maps $h:U \rightarrow Y$ cover $Y$, one has internally in $\Tcal$: 
``$\forall y \in Y, \exists h \in H$ such that $y$ is in the image of $h$", where $H$ denote the external set of such map $h:U\rightarrow Y$ which are finite and of bounded cardinal in $\Tcal_{/Y}$. In particular the joint image of all the $\phi_h$ for $h \in H$ contains all the generators of $l^{2}(Y)$ and hence spam a dense subspace of $l^{2}(Y)$, which concludes the proof by proposition \ref{QcontaininHcal}.
}

}

\blockn{We denote by $\Hcal^{red}(\topos)$ the full subcategory of $\Hcal(\topos)$ of objects which are weakly contained in $l^{2}(X)$ for some separating object $X$. Objects of $\Hcal^{red}(\topos)$ are said to be \emph{square-integrable}\footnote{because when $\topos$ is the topos of $G$-sets for some discrete group $G$, then $\Hcal(\Tcal)$ is the category of unitary representations of $G$ while $\Hcal^{red}(\topos)$ is precisely the category of square integrable representations of $G$}.

Results of \cite{henry2014measure} (especially section 7) suggest that the square integrable Hilbert spaces of $\Tcal$ are the one that are clearly related to the geometry of $\topos$. 

Because of proposition \ref{l2weaklycontinl2}, for any separating bound $X$ of $\topos$ the Hilbert space $l^{2}(X)$ is a generator of $\Hcal^{red}(\topos)$, hence extending proposition 7.6 of \cite{wstarcat} to monotone complete $C^{*}$-category (or restricting ourselves to topos which are integrable in the sense of \cite[section 3]{henry2014measure}) gives us that when $\Tcal$ is a boolean and locally separated topos, $\Hcal^{red}(\topos)$ is equivalent to the category of reflexive Hilbert modules over $End(l^{2}(X))$. For this reason we can call this algebra ``the" (or ``a") \emph{reduced}\footnote{For example, if $\Tcal$ is the topos of $G$-sets for $G$ a discrete group one obtains the usual (reduced) von Neumann algebra of the group this way.} algebra of $\topos$ (it is unique up to Morita equivalence).
}

\blockn{Finally, The fact that it is possible to define internally the tensor product of two Hilbert spaces yields a symmetric monoidal structure on $\Hcal(\topos)$. Moreover as $l^{2}(X) \otimes l^{2}(Y) \simeq l^{2}(X \times Y)$ one can see that $\Hcal^{red}(\topos)$ is stable by tensor product\footnote{in \ref{tensorbySI} we will actually prove the stronger result that the tensor product of an arbitrary Hilbert space with a square integrable Hilbert space is square integrable}. Hence, $\Hcal^{red}(\topos)$ is also endowed with a symmetric monoidal structure, but without a unit object (unless $\Tcal$ is separated). }

\section{Statement of the main theorems}
\label{secStatement}
\block{\Def{Let $C$ be a monotone complete $C^{*}$-category, $\Ecal$ any topos (not necessary boolean) a representation of $C$ in $\Ecal$ is a $*$-functor $\rho$ from $C$ to $\Hcal(\Ecal)$. It is said to be normal if for any supremum $a = \sup a_i$ of a bounded directed net of positive operator in $C$, $\rho(a_i)$ converge internally in the weak operator topology to $\rho(a)$.}

A representation of $C$ in $\Ecal$ is also the same as a representation of $p^{*} C$ in the category of Hilbert space internally in $\Ecal$ (where $p$ is the geometric morphism from $\Ecal$ to the point). Also if $\Ecal$ is boolean, then $\Hcal(\Ecal)$ is monotone complete and a representation $\rho$ is normal in the sense of this definition if and only if it is normal as a $C*$-functor between $C^{*}$-category. }

\block{When $C$ has additional structure (for example is monoidal) we will by default assume that the representation $\rho$ preserve these structures, for exemple:

\label{DefRepUnred}\Def{If $\Tcal$ is a boolean topos and $\Ecal$ an arbitrary topos, a representation of $\Hcal(\Tcal)$ in $\Ecal$ is a normal representation of the monotone complete $C^{*}$-category $\Hcal(\Tcal)$ in $\Ecal$ such that the underlying $*$-functor is symmetric and monoidal. }}

\block{\Def{We will say that a Hilbert space $H \in \Hcal(\Ecal)$ is inhabited if one has $\exists s \in H$ such that $\Vert s \Vert >0$ internally in $\Ecal$, or equivalently, if $\exists s \in H, \Vert s \Vert =1$ holds internally in $\Ecal$.}
}

\block{\label{DefNondegenUnred}\Def{If $\Tcal$ is a boolean topos and $\Ecal$ an arbitrary topos, a representation $\rho$ of $\Hcal(\Tcal)$ is said to be non-degenerate if for any inhabited object $H$ of $\Hcal(\Tcal)$ the Hilbert space $\rho(H)$ is inhabited in $\Ecal$.}

It is important to notice that detecting whether a representation $\rho$ of $\Hcal(\Tcal)$ is non-degenerate or not can be done completely from the (monoidal) category $\Hcal(\Tcal)$ without knowing the topos $\Tcal$. Indeed, an object $H$ of $\Hcal(\Tcal)$ is inhabited in $\Tcal$ if and only if the functor $H \otimes \_$ is faithful. 

As an example of a ``degenerate'' representation, one can consider $\Tcal$ the topos of $G$-sets for some infinite discrete groupe $G$, then $\Hcal(\Tcal)$ is the category of unitary representation of $G$. Let $\rho$ be the representation of $\Hcal(\Tcal)$ into $\Tcal$ defined by:

\[\rho(H) = \overline{ \{ h \in H | h \text{ belongs to a finite dimensional sub-representation } \} } \]

One has $\rho(l^{2}(G))= 0$ so it cannot be non-degenrate. One easily checks that it is a symmetric normal $*$-functor, and it is monoidal because of the following observation:

\Lem{Let $R,R' \in \Hcal(\Tcal)$ be two representations of $G$. Assume that $R \otimes R'$ contains a (non trivial) finite dimensional sub-representation, then both $R$ and $R'$ contains a non trivial finite dimensional sub-representation. }

\Dem{ Let $K \subset R \otimes R'$ a non-trivial finite dimensional sub-representation. let $K^{*}$ the dual of the representation $K$. Because $K$ is finite dimensional and non-trivial $K \otimes K^{*}$ contains a non-zero invariant vector. In particular $R \otimes R' \otimes K^{*}$ contains a non-zero invariant vector which corresponds to a non zero $G$-linear Hilbert-Schmidt $f$ from $R$ to $R'^{*}\otimes K$, $f^{*}f$ is hence a non-zero compact self-adjoint $G$-linear automorphism of $R$ which is hence going to have some non-trivial finite dimensional $G$-stable eigenspaces. This concludes the proof of the lemma.}

}

\blockn{There is in particular a representation of $\Hcal(\Tcal)$ in $\Tcal$ itself called the \emph{tautological} representation and given by the identity functor. This representation is non-degenerate. Moreover if $\rho$ is a (non-degenerate) representation of  $\Hcal(\Tcal)$ in a topos $\Ecal$ and $f:\Ecal' \rightarrow \Ecal$ is any geometric morphism then $f^{\sharp}\rho$, given by composing the functor $\rho$ by the functor $f^{\sharp} : \Hcal(\Ecal) \rightarrow \Hcal(\Ecal')$ is a (non-degenerate) representation of $\Hcal(\Tcal)$ in $\Ecal'$. In particular, any geometric morphism from $f:\Ecal \rightarrow \Tcal$ induce a non-degenerate representation $f^{\sharp}$ of $\Hcal(\Tcal)$ in $\Ecal$.}

\block{\label{mainThunred}\Th{Let $\Tcal$ be a boolean locally separated topos, then $\Tcal$ is the classifying topos for non-degenerate representations of $\Hcal(\Tcal)$ and the tautological representation is the universal non-degenerate representation. More precisely, for any topos $\Ecal$ there is an equivalence of category from the category of geometric morphism $\hom(\Ecal,\Tcal)$ to the category\footnote{See \ref{Defmorphofrep} below.} of non-degenerate representations of $\Hcal(\Tcal)$ in $\Ecal$ which associate to any geometric morphism $f$ the representation $f^{\sharp}$.  }

In particular, the topos $\Tcal$ is uniquely determined (up to unique isomorphism) from the symmetric monoidal $C^{*}$-category $\Hcal(\Tcal)$ of Hilbert bundles over $\Tcal$.

}

\block{\label{Defmorphofrep}This theorem, as it is stated, does not completely make sense yet because we did not say what are the morphisms of representations of $\Hcal(\Tcal)$ in $\Ecal$. It appears that the good notion of morphisms is the following:

\Def{If $\rho_1$ and $\rho_2$ are two representations of $\Hcal(\Tcal)$ in $\Ecal$ a morphism $v$ from $\rho_1$ to $\rho_2$ is a collection of isometric inclusions $v_X:\rho_1(X) \rightarrow \rho_2(X)$ for each object $X$ of $\Hcal(\Tcal)$ which is natural in $X$ and such that, up to the structural isomorphisms, $v_{\C_{\Tcal}}$ is the identify of $\C_{\Ecal}$ and $v_X \otimes v_Y$ is $v_{X \otimes Y}$. }

Because we do not assume that the $v_X$ have adjoints they are not morphisms in $\Hcal(\Ecal)$ and hence it would not make sense strictly speaking to says they form a symmetric monoidal natural transformation, but this is essentially what this definition means.

\bigskip

It can be proved directly, using the fact that $\Hcal(\Tcal)$ as a generator (as proved in \cite[3.6]{henry2014measure}), that morphisms between two representations $\rho_1$ and $\rho_2$ actually form a set and not a proper class, but we do not need to know that and one can instead obtain this result as a corollary of the theorem.
}

\block{There is also a form of this theorem for the category $\Hcal^{red}(\Tcal)$ instead of $\Hcal(\Tcal)$, which we will use as an intermediate step in the proof of theorem \ref{mainThunred}. The problem to state it directly is that $\Hcal^{red}(\Tcal)$ is not exactly a monoidal category because it does not have a unit object. In particular, one cannot ask monoidal functors to preserve the unit object, but it appears that the condition ``non-degenerate" actually completely replace the preservation of the unit object, and moreover, in this case, the definition of non-degenerate can be weakened:

\label{DefRepReduced}\Def{A non-degenerate representation $\rho$ of $\Hcal^{red}(\topos)$ in $\topos[1]$ is a normal representation of the monotone complete $C^{*}$-category $\Hcal^{red}(\topos)$ in $\Ecal$ such that:

\begin{itemize}

\item $\rho$ satisfies all the axioms of the definition of a symmetric monoidal functor not involving the unit object.

\item There exists an object $H$ of $\Hcal^{red}(\topos)$ such that $\rho(H)$ is inhabited.

\end{itemize}
}

The morphisms of such representations are defined exactly as in \ref{Defmorphofrep}.}

\blockn{Similarly to \ref{mainThunred}, there is a tautological representation of $\Hcal^{red}(\Tcal)$ in $\Tcal$ and it is possible to pullback any representation of $\Hcal^{red}(\Tcal)$ along a geometric morphism }

\block{\label{mainThRed}\Th{The tautological representation of $\Hcal^{red}(\topos)$ in $\topos$ is the universal non-degenerate representation of $\Hcal^{red}(\topos)$. I.e., for any topos $\topos[1]$ the association $p \rightarrow p^{*}$ induce an equivalence of categories between geometric morphisms from $\topos[1]$ to $\topos$ and non-degenerate representations of $\Hcal^{red}(\topos)$ in $\topos[1]$.}

}

\blockn{In the rest of the article, all the representations considered will always be implicitly assumed to be non-degenerate.}

\blockn{These two theorems together, suggest that it should be possible to reconstruct $\Hcal(\Tcal)$ directly from $\Hcal^{red}(\Tcal)$. This indeed seems to be the case, here is a sketches of proof:

We will see in \ref{tensorbySI}, that any $H$ in $\Hcal(\Tcal)$ induce by tensorisation an endofunctor:

 \[ M_H:= ( H \otimes \_ )  : \Hcal^{red}(\Tcal) \rightarrow \Hcal^{red}(\Tcal) \]
 
which satisfies a ``multiplier" condition of the form $M_H(A\otimes B) \simeq M_H(A) \otimes B$ (this isomorphism being functorial and satisfying some coherence conditions). Conversely, if $M$ is an endofunctor of $\Hcal^{red}(\Tcal)$ satisfying the same condition than $M_H$, then for any separating object $X$ of $\Hcal$, $M_H(l^{2}(X))$ is going to be a $l^{2}(X)$-module in the sense of definition \ref{Defl2module}, hence, by proposition \ref{l2Xmodeqtoslice}, is of the form $\bigoplus_{x \in X} H_x$ for a Hilbert space $(H_x)_{x \in X}$ in $\Tcal_{/X}$, but using the coherence condition on the isomorphisms $M(A\otimes B) \simeq M(A) \otimes B$ one should be able to prove that all the $H_x$ are canonically isomorphic and hence (if $X$ is inhabited) that $M(l^{2}(X))$ if of the form $H \otimes l^{2}(X)$. One can then conclude the $M$ is exactly the tensorization by $H$ by assuming that $X$ is a separating bound and hence that $l^{2}(X)$ is a generator of $\Hcal^{red}(\Tcal)$. One also has to check that endomorphisms of $l^{2}(X)$ acts as they should on $M(l^{2}(X)) \simeq H \otimes l^{2}(X)$ but this will follow by the same ``matrix elements" argument as in the proof of \ref{proofcompatibilityRED}.

\bigskip

Hence, $\Hcal(\Tcal)$ should be, in some sense, the ``category of multiplier" of the (non-unital) monoidal category $\Hcal^{red}(\Tcal)$.

\bigskip

We decided not to include a precise form of this last result in the present version of this paper, because its proof, and in fact even its proper formulation, require some work on the precise coherence conditions required on such a ``multiplier" that seems to be out of the scope of the present paper, and we are not sure that such a statement has, despite its elegance, any interesting applications.

}

\section{From representations of $\Hcal^{red}$ to geometric morphisms}
\label{secProofRed}
\renewcommand{\thesubsubsection}{\arabic{section}.\arabic{subsection}.\arabic{subsubsection}}

\blockn{In this section we consider an arbitrary representation $\rho$ of $\Hcal^{red}(\topos)$ in a topos $\topos[1]$, and we will prove that it induces a geometric morphism from $\topos[1]$ to $\topos$. In all this section we will work internally in $\topos[1]$. Hence objects of $\topos[1]$ will be called sets and objects of $\Hcal(\topos[1])$ will simply be called Hilbert spaces, the price of this being to avoid the use of the law of excluded middle and of the axiom of choice. }

\subsection{Construction of the geometric morphism on separating objects}
\label{GeomorphonObjec}
\block{In this subsection, we fix a separating object $X \in \topos$, and we will show that $\rho(l^{2}(X))$ is of the form $l^{2}(Y)$ for a well determined decidable set $Y$. We will denote by $H$ the space $\rho(l^{2}(X))$. This space is endowed with an operator $ \Delta : H \rightarrow H \otimes H$ which comes from the operator in $\Hcal(\topos)$:

\[ \begin{array}{c c c c}
\Delta : & l^{2}(X) &\rightarrow &l^{2}(X) \otimes l^{2}(X) \\
 & e_x & \mapsto & e_x \otimes e_x \end{array} \]

$\Delta$ is co-commutative and co-associative and its adjoint $\Delta^{*} : H \otimes H \rightarrow H$ defines a multiplication on $H$ denoted $ a * b :=\Delta^{*}(a \otimes b)$, which is commutative and associative. Moreover $\Delta^{*} \Delta = Id_H$ hence $\Vert a * b \Vert \leqslant \Vert a \Vert \Vert b \Vert$.

In $\topos$, the operator $\Delta^{*}$ is defined by $\Delta^{*}(e_x \otimes e_y)=e_x$ is $x=y$ and $0$ if $x \neq y$, or, tu put it another way, it is just the point-wise multiplication of $X$-indexed sequences.

In $\Ecal$, for $ h \in H$, we denote by $m_h$ the linear map on $H$ defined by $m_h(u)=h * u$. we denote by $\Vert h \Vert_{\infty} = \Vert m_h \Vert$, and the usual norm of $H$ will be denoted by $\Vert h \Vert_2$. One has $\Vert h \Vert_{\infty} \leqslant \Vert h \Vert_{2}$ because of the previous bound on the norm of a product.

Finally, let $C$ be the closure in the space of bounded linear map from $H$ to $H$ of the of algebra formed by the $m_h$ for $h \in H$. It is a commutative algebra. the term ``closure" is taken here in the sense that:

\[ C := \{ f :H \rightarrow H | \forall \epsilon > 0, \exists h \in H, \Vert f - m_h \Vert < \epsilon \} \]

In fact, it should be called ``fiberwise closure" or ``weak closure" because the resulting set is not closed in the sense the complementary of an open, see \cite{henry2014localic} for more details about these notions.

}

\block{\label{X0}Before going further we need to make a few external constructions that will also be usefull in the next subsetions.

Because there exists an object $H$ of $\Hcal(\Tcal)$ such that $\rho(H)$ is inhabited, and that the $l^{2}(X)$ for $X$ separating are generators of $\Hcal(\Tcal)$ there exists a separating object $X_0 \in \Tcal$ such that $\rho(l^{2}(X_0))$ is also inhabited. We fix such an object $X_0$.

A sub-object $U \subset X_0 \times X$ is said to be of degree smaller than $n$ over $X_0$ if internally in $\topos$ for all $x \in X_0$ there is at most $n$ distinct elements $x' \in X$ such that $(x,x') \in U$, i.e. if $U$ is of cardinal smaller than $n$ as an object of $\Tcal_{/X_0}$. It is said to be of finite degree if it is of degree smaller than $n$ for some \emph{external} natural number $n$, i.e. if, as an object of $\Tcal_{/X_0}$ its cardinal is bounded.

Theorem \ref{SepImpFinitness} applied to $\Tcal_{/X_0}$ imply that sub-objects of finite degree over $X_0$ cover $X_0 \times X$. We denote by $F$ the (external) ordered set of sub-objects of $X_0 \times X$ of finite degree. As, if $U$ and $U'$ are of finite degree then $U \vee U'$ is also of finite degree, $F$ is a directed pre-ordered set. }

\block{\label{VUPU}For $U \in F$ let (in $\topos$):

\[ \begin{array}{c c c c}
V_U :& l^{2}(X_0) &\rightarrow & l^{2}(X_0) \otimes l^{2}(X) \\
& e_x &\mapsto & \displaystyle \sum_{x' \text{ s.t. } (x,x')\in U} e_x \otimes e_{x'}
\end{array}
\]

\[ \begin{array}{c c c c}
P_U :& l^{2}(X_0) \otimes l^{2}(X) &\rightarrow & l^{2}(X_0) \otimes l^{2}(X) \\
& e_x \otimes e_{x'} & \mapsto &\displaystyle  (\Char_{(x,x') \in U } ) e_x \otimes e_{x'}
\end{array}
\] 

We also denote by $V_U$ and $P_U$ their images by $\rho$, and we denote $H_0$ the space $\rho(l^{2}(X_0))$.

One easily check that:

\begin{itemize}

\item $V_U$ has an adjoint given by $V_U^{*}(e_x \otimes e_y) =e_x \Char_{(x,y) \in U}$

\item $V_U$ is well defined and bounded (if $U$ has degree $n$, then $V_U$ has norm smaller than $\sqrt{n}$).

\item $P_U$ is a $F$-indexed net of projections whose supremum is the identity of $l^{2}(X_0) \otimes l^{2}(X)$, and as $\rho$ is normal it is also the case in $\topos[1]$ that the supremum of the $P_U$ is the identify of $H_0 \otimes H$.

\item $(Id_{H_0} \otimes \Delta^{*}) \circ (V_U \otimes Id_H )= P_U$. Indeed this is easily checked internally in $\topos$ on elements of the form $e_x \otimes e_{x'}$, extended by linearity and continuity and transported to $\topos[1]$ by $\rho$.

\end{itemize}
}

\block{\label{m_injectif}The role of these operators $V_U$ and $P_U$ is essentially to provide ``locally" an approximate unit for the product $*$. Indeed, let us fix an element $h_0 \in H_0$ such that $\Vert h_0\Vert =1$, which is always possible by assumption on $X_0$, then one can define $1_U := (h_0^{*} \otimes Id_H) V_U(h_0)$, and, because of the formula relating $V_U$ and $P_U$ one has after a short computation that for any $h \ \in H$:

\[ h*1_U =  (h_0^{*} \otimes Id_{H}) (P_U(h_0 \otimes h)), \] 

which (because $P_U$ converge strongly to the identity) converge in norm to $h$. In particular:

\Lem{In $\Ecal$, the association $h \mapsto m_h$ from $H$ to $C$ is injective.}

Indeed, using (internally) our approximate unit $1_U$ one gets that:

\[ h= \lim_{U \in F} h * 1_U =\lim_{U \in F} m_h(1_U). \]

}

\block{We can now prove:

\Lem{The algebra $C$ generated by the $m_h$ is a $C^{*}$-algebra.}

The term $C^{*}$-algebra is taken here in the same sense as, for example, in \cite{banaschewski2000spectral}. I.e. it is complete in the sense of Cauchy approximation or Cauchy filters, and the ``norm" is describe either by the data of the ``rational ball" $B_q$ corresponding morally to the set of $x$ such that $\Vert x \Vert x < q$ of equivalently as a function $\Vert \_ \Vert$ with value into the object of upper semi-continuous real number.

\Dem{ All we have to do is to prove that every element of $C$ have an adjoint which belongs to $C$. As the adjunction preserve the norm and that the space of bounded linear map is complete (in the norm topology) it is enough to prove it for a dense family of elements of $C$.

Let $h$ be an arbitrary element of $H$.
\bigskip

We define: 
 
\[h_U = h * 1_U = (h_0^{*} \otimes Id_{H}) (P_U(h_0 \otimes h)) \in H \]

As mentioned earlier, $(h_U)_{U \in F}$ converges in norm to $h$.

We also define:

\[h'_U = ( h_0^{*} \otimes h^{*} \otimes Id_H) (Id_{H_0} \otimes \Delta)  (V_U)(h_0)  \]

We will prove that for all $a,b \in H$ one has:

\[ \scal{h'_U * a}{b} = \scal{a}{h_U*b} \]

this will show that $m_{h_U}$ has an adjoint $m_{h'_U}$ which belong to $C$. As for any $h\in H$, $m_h$ is approximated by the $m_{h_U}$ (indeed, $\Vert m_h \Vert = \Vert h \Vert_{\infty} \leqslant \Vert h \Vert _2$ ) the $m_{h_U}$ for $h \in H$ and $U \in F$ are dense in $C$, hence this will conclude the proof.

In order to do so, consider the following two operators $W_1$ and $W_2$ from $H_0 \otimes H \otimes H$ to $H_0 \otimes H$ defined by:

\[ W_1 := (Id_{H_0} \otimes \Delta^{*}) \circ (Id_{H_0} \otimes \Delta^{*} \otimes Id_{H}) \circ (V_U \otimes Id_H \otimes Id_H) \]

\[W_2:= (V_U^{*} \otimes Id_H) \circ (Id_{H_0} \otimes \Delta^{*} \otimes id_H) \circ (Id_{H_0} \otimes Id_H \otimes \Delta) \]

Two short computations show that:
\[ \scal{a}{h_U * b} = \scal{h_0 \otimes a }{W_1(h_0 \otimes h \otimes b)} \]
\[ \scal{h'_U * a}{b} = \scal{h_0 \otimes a }{W_2(h_0 \otimes h \otimes b)} \]

Moreover, both $W_1$ and $W_2$ are image by $\rho$ of operators also denoted $W_1$ and $W_2$ defined by the same formulas in $\Tcal$ (replacing $H_0$ and $H$ by $l^{2}(X_0)$ and $l^{2}(X)$) and one can easily check in $\Tcal$ that $W_1 = W_2$ indeed, if $x \in X_0, y,z \in X$ a short computation shows that:

\[ W_1 (e_x \otimes e_y \otimes e_z) = \mathbb{I}_{(y=z)} \mathbb{I}_{((x,y) \in U)}e_x \otimes e_y =  W_2 (e_x \otimes e_y \otimes e_z)  \]

which proves that $W_1 = W_2$ in $\Tcal$ and hence also in $\Ecal$ and hence that $\scal{a}{h_U * b} = \scal{h'_U * a}{b}$ which concludes the proof as mentioned above.

}

}

\block{\label{geomofC}One of the key observation is that this algebra $C$ is ``geometrically" attached to the representation $\rho$, in order to make this clear let us denote $C_{\rho}$ instead of $C$ in what follows.

\Lem{Let $p:\Ecal' \rightarrow \Ecal$ be any topos above $\Ecal$. One has a representation $p^{\sharp}\rho$ of $\Hcal(\Tcal)$ in $\Ecal'$, and a canonical isomorphism:

\[ C_{p^{\sharp} \rho} \simeq p^{\sharp} C_{\rho}. \]

}

\Dem{It is clear that $p^{\sharp}\rho$, obtained as the composition of $\rho :\Hcal^{red}(\Tcal) \rightarrow \Hcal(\Ecal)$ with $p^{\sharp} : \Hcal(\Ecal) \rightarrow \Hcal(\Ecal')$, is a representation in the sense of definition \ref{DefRepReduced}. Moreover, the $C^{*}$-algebra $p^{\sharp}(C_{\rho})$ is naturally acting on $p^{\sharp}(H) = p^{\sharp}\rho (l^{2}(X))$ simply because $C_{\rho}$ is acting on $H$. Moreover if $h \in H$, the operator $p^{\sharp}(m_h)$ on $p^{\sharp}(H)$ is equal to the operator $m_{p^{*}(h)}$ (this is clear from the definition of $m_h$). 
Hence the image of $p^{\sharp}(C_{\rho})$ in $End(p^{\sharp}(H))$ is generated by operators of the form $m_{h}$ for $h \in p^{*}(H) \subset p^{\sharp}(H)$ but as $p^{*}(H)$ is dense in $p^{\sharp}(H)$ this proves that the image of $p^{\sharp}(C_{\rho})$ in the algebra of endomorphism of $p^{\sharp}(H)$ is indeed the algebra generated by the $m_h$ for $ h \in p^{\sharp}(H)$, i.e. the algebra $C_{p^{\sharp}\rho}$, giving a canonical surjection from $p^{\sharp} C_{\rho}$ to $C_{p^{\sharp}\rho}$.  But this surjection is also isometric and hence is an isomorphism because of the relation:

\[ \Vert m_h \Vert = \lim_{n \rightarrow \infty} \Vert h * h* \dots *h \Vert ^{\frac{1}{n}}\]  

which is proved in the proof of lemma \ref{LemChascontnorm} below and clearly preserved by $p^{\sharp}$.
}

In particular, if one compute $C$ in the tautological representation of $\Hcal^{red}(\Tcal)$ in $\Tcal$ it is exactly $\Ccal_0(X)$, hence if our representation $\rho$ corresponds to a pullback along some geometric morphism $f$ the previous lemma imply that $C$ will be $f^{\sharp}(\Ccal_0(X)) = \Ccal_0( f^{*}(X))$ and hence using the non-unital Gelfand duality proved in \cite{henry2014nonunital} one has that $f^{*}(X)$ can be reconstructed out of $\rho$ as the spectrum of $C$. In order to show that an arbitrary $\rho$ fo comes from a geometric morphism one needs to prove that the spectrum of $C$ is a discrete decidable locale, and we will do that by showing that $\spec C$ is locally positive (i.e. that $\spec C \rightarrow 1$ is an open map) and that the diagonal map $\spec C \rightarrow \spec C \times \spec C$ is both open (hence by \cite[C3.1.15]{sketches} that $\spec C$ is discrete) and closed (hence that $\spec C$ is decidable).

}

\block{
\label{LemChascontnorm}\Lem{The norm of any element of $C$ is a continuous real number and $\spec C$ is locally positive (or ``open", or ``overt").}

\Dem{Let $h \in H$ such that $\Vert h \Vert_2>0$, as $\lim_U m_h (1_U)= h$ this proves that $\Vert m_h \Vert >0$ in the sense that there exists a rational number $q$ such that $0<q \leqslant \Vert m_h\Vert$. Let $h^{(n)}$ denotes the $n$-th power of $h$ for the product $*$. As $m_h$ is an element of a commutative $C^{*}$-algebra one has $\Vert m_h^{n} \Vert = \Vert m_h \Vert^{n}$ hence $\Vert h^{(n)}\Vert_{\infty} = \Vert h \Vert_{\infty}^{n}$.

Moreover:

\[ \Vert h^{(n)} \Vert_2 = \Vert m_h h^{(n-1)} \Vert_2 \leqslant \Vert m_h \Vert \Vert h^{(n-1)} \Vert_2 \leqslant \Vert h \Vert_{\infty} \Vert h^{(n-1)} \Vert_2 \]

Hence by induction on $n$: 
\[ \Vert h^{(n+1)} \Vert_2 \leqslant \Vert h \Vert_{\infty}^{n} \Vert h \Vert_2 \]

As $\Vert . \Vert_{\infty} \leqslant \Vert . \Vert_{2}$ one has:

\[ q\Vert h \Vert_{\infty}^{n} \leqslant \Vert h \Vert_{\infty} ^{n+1} \leqslant \Vert h^{(n+1)} \Vert_2 \leqslant \Vert h \Vert_{\infty}^{n} \Vert h \Vert_2 \]

For some rational number $q$ such that $0<q<\Vert h \Vert_{\infty}$. Taking the $n$-th root one obtains that:

\[ q^{\frac{1}{n}} \Vert h \Vert_{\infty} \leqslant \Vert h^{(n+1)} \Vert_2^{\frac{1}{n}} \leqslant \Vert h \Vert_{\infty} \Vert h \Vert_2^{\frac{1}{n}} \leqslant \Vert h \Vert_{\infty} (q')^{\frac{1}{n}} \]

And equivalently that:

\[\frac{1}{q'^{1/n}} \Vert h^{(n+1)} \Vert_2^{\frac{1}{n}}  \leqslant  \Vert h \Vert_{\infty} \leqslant \frac{1}{q^{1/n}} \Vert h^{(n+1)} \Vert_2^{\frac{1}{n}}\]

As $\Vert h^{(n+1)} \Vert_2^{\frac{1}{n}}$ is bounded, the difference between the upper bound and the lower bound can be made arbitrarily small, hence $\Vert h \Vert_{\infty}$ can be approximated arbitrarily close by continuous real number (and hence also by rational number), which proves that $\Vert h \Vert_{\infty}$ is a continuous real number. One also gets the identity:

\[ \lim_{n \rightarrow \infty} \Vert h^{(n+1)} \Vert_2^{\frac{1}{n}} = \Vert h \Vert_{\infty} \]

which was required in the proof of the previous lemma.

As $0$ and the $h$ such that $\Vert h \Vert_2 >0$ are dense in $H$ (for any $h \in H$ and any $\epsilon>0$ either $h$ is of positive norm or it is of norm smaller than $\epsilon$ and hence is approximated by $0$) and as the $m_h$ for $h \in H$ are dense in $C$, one has obtain a dense familly of elements of $C$ of continuous norm hence every element of $C$ has a continuous norm and by \cite[5.2]{henry2014nonunital} one can conclude that $\spec C$ is locally positive.
}}

\block{\label{const_e_chi}If one works internally in $\spec C$, one still has a representation $\rho$ of $\Hcal(\Tcal)$ (by pulling it back from the one in $\Ecal$) and by \ref{geomofC} the $C^{*}$-algebra of endomorphisms of $\rho(l^{2}(X))$ generated by this representation is isomorphic to the pullback of $C$, hence (by \cite[proposition 4.4]{henry2014nonunital} ) we have at our disposal a character $\chi : C \rightarrow \C$ satisfying $\exists c \in C$ such that $|\chi(c)|>0$. We will now examine what can be done with such a character, without necessary assuming that we are working internally in $\spec C$.

\bigskip

Assume one has a character $\chi$ of the $C^{*}$ algebra $C$ satisfying $\exists c \in C, |\chi(c)|>0$. The map $h \rightarrow \chi(m_h)$ is a bounded linear form on $H$, also denoted $\chi$. Also there exists a $g \in H$ such that $\chi(g)=1$.

Let $e_\chi$ be defined as $(\chi \otimes Id_H)(\Delta(g))$.

Then if $h$ is another element of $H$ one has:

\[ \scal{h}{e_\chi} = (\chi \otimes h^{*})(\Delta(g)) = \chi( (Id_H \otimes h^{*})(\Delta(g))) \]

But, for any $w \in H$ one has:

\[ \scal{w}{(Id_H \otimes h^{*})(\Delta(g))} = (w^{*} \otimes h^{*})(\Delta(g)) = \scal{\Delta^{*}(w \otimes h)}{g} = \scal{m_h w}{g} \]

hence $(Id_H \otimes h^{*})(\Delta(g)) = (m_h)^{*}(g)$ and as $\chi$ is a character of $C$ one obtain that $\chi( (Id_H \otimes h^{*})(\Delta(g))) = \overline{\chi(h)}\chi(g)=\overline{\chi(h)}$. Finally:

\[ \scal{e_\chi}{h} = \chi(h) \]

In particular $e_{\chi}$ does not depend on $g$, but only on $\chi$. This proves:

\Lem{If $\chi$ is a character of $C$ then there exists a unique element $e_\chi \in H$ such that for all $h \in H$ one has $\scal{e_\chi}{h} = \chi(m_h)$.
}

}

\block{\label{e_chi_prop}The vector $e_{\chi}$ attached to a character $\chi$ of $C$ satisfies additional properties:

\Lem{\begin{enumerate}
\item $ \Delta(e_{\chi}) = e_{\chi} \otimes e_{\chi} $
\item $ e_{\chi} * e_{\chi} = e_{\chi} $
\item $ e_{\chi}*h=\chi(h) e_{\chi} $
\item $ \Vert e_{\chi} \Vert_2 = \Vert e_{\chi} \Vert_{\infty} = 1 $
\end{enumerate}}

\Dem{
\begin{enumerate}
\item Let $a,b \in H$, then:

\[\scal{\Delta(e_{\chi})}{a \otimes b} = \scal{e_{\chi}}{a * b} = \chi(a * b)=\chi(a) \chi(b)= \scal{e_{\chi} \otimes e_{\chi}}{a \otimes b} \] which proves the results.

\item One has $e_{\chi} * e_{\chi} = \Delta^{*}(e_{\chi} \otimes e_{\chi}) = \Delta^{*}\Delta(e_{\chi}) = e_{\chi}$ because $\Delta^{*}\Delta =Id$.

\item For any $h,g \in H$ one has: 
\[\scal{e_{\chi}*h}{g} = \scal{m_h (e_\chi)}{g} = \scal{e_\chi}{(m_h)^{*} g } =\chi((m_h)^{*} m_g)=\overline{\chi(h)}\chi(g).\]

In particular $\scal{e_{\chi}*h}{g} = \scal{\chi(h)e_{\chi}}{g}$ which concludes the proof.

\item $\scal{e_\chi}{e_\chi} = \chi(e_{\chi})$. But as $e_{\chi}$ is a projector for the product $*$, the value of $\chi(e_{\chi})$ is either $0$ or $1$. The norm $\Vert e_{\chi} \Vert_2$ of $e_{\chi}$ is hence either $0$ or $1$, but it cannot be $0$ (because $\chi$ is non-zero) hence it is one. This implies that $\Vert e_{\chi} \Vert_{\infty} \leqslant 1$ but has $m_{e_{\chi}} e_{\chi} = e_{\chi}$ one has $\Vert e_{\chi} \Vert_{\infty} =1$.

\end{enumerate}

}

}

\block{\label{specC_atomic}\Lem{The diagonal of $\spec C$ is both open and closed.}

\Dem{Working internally in $\spec C \times \spec C$ one has two characters $\chi_1$ and $\chi_2$ which give rise to elements $e_{\chi_1}$ and $e_{\chi_2}$ of (the completion of the pullback of) $H$ as in \ref{const_e_chi}.

\[ \scal{e_{\chi_1}}{e_{\chi_2}} = \chi_1(e_{\chi_2}) \]

Hence, as $e_{\chi_2}$ is a projector, $\scal{e_{\chi_1}}{e_{\chi_2}}$ is either $0$ or $1$. If it is $0$ then $\chi_1 \neq \chi_2$ but if it is one, as $\Vert e_{\chi_1}\Vert_2 =\Vert e_{\chi_2}\Vert_2 =1$ this imply that $\chi_1=\chi_2$. Hence, internally in $\spec C \times \spec C$ one has $\chi_1=\chi_2$ or $\chi_1 \neq \chi_2$ which proves that the diagonal of $\spec C \times \spec C$ is both open and closed.
}

}

\block{\Prop{There exists a unique (up to unique isomorphism) decidable set $X_{\rho}$ such that $H=l^{2}(X_{\rho})$ and $\Delta$ is the map defined by $\Delta(e_x)=e_x \otimes e_x$ for all $x \in X_{\rho}$. }

\Dem{
 The uniqueness of $X_{\rho}$ is clear: any (decidable) object $X$ can be reconstructed out of $l^{2}(X)$ and $\Delta$ as the set of $v \in l^{2}(X)$ such that $\Vert v \Vert = 1$ and $\Delta(v) = v \otimes v$.

As $\spec C$ is open and its diagonal map is open it is a discrete locale by \cite[C3.1.15]{sketches}, i.e. $\spec C = X_{\rho}$ for some set $X_{\rho}$. As the diagonal of $\spec C$ is also closed, $ X_{\rho}$ is decidable. Each $x \in X_{\rho}$ corresponds to a character of $C$ and hence there is a corresponding vector $e_x \in H$ such that $\Delta(e_x)=e_x \otimes e_x$. Moreover the proof of \ref{specC_atomic} show that $(e_x)$ is an orthonormal family in $H$, hence $l^{2}(X_{\rho})$ embeds isometrically in $H$ in a way compatible to $\Delta$.

Let $x \in X_{\rho}$. As an element of $C \simeq \Ccal_0(X_{\rho})$, $m_{e_x}$ is the characteristic function of $\{x\}$ : indeed $\chi_y(m_{e_x}) = \scal{e_y}{e_x}=\delta_{x,y}$. Hence the characteristic function  of $\{ x \}$ acts on $H$ by $m_{e_x} h = \scal{e_x}{h} e_x$ and has its image included in $l^{2}(X_{\rho})$. But any element of $C = \Ccal_0(X_{\rho})$ is a norm limit of linear combination of such characteristic function, hence the action of $C$ has its image included in $l^{2}(X_{\rho})$. But has this action contains $m_{1_U}$ which converge to the identity, this implies that $l^{2}(X_{\rho}) = H$ (in a way compatible to $\Delta$ ).
}
}

\subsection{Functoriality on separating object}
\label{secFunctoOnObj}
\block{\label{DefGammaf}Let $X$ and $Y$ be two decidable objects of any topos, and let $f:X \rightarrow Y$ be a map. Then the map $f$ induces an operator:

\[\begin{array}{c c c c}
\Gamma_f :& l^{2}(X) & \rightarrow & l^{2}(X) \otimes l^{2}(Y) \\
& e_x & \mapsto & e_x \otimes e_{f(x)}
\end{array} \]

indeed, it has an adjoint $\Gamma_f^{*}$ which send $e_x \otimes e_y$ to $e_x$ if $y=f(x)$ and to zero otherwise (which make sense constructively because $Y$ is decidable).

}

\block{\label{functorial} If $X$ and $Y$ are two objects of $\Tcal$, we already know (from the previous section) that there exists two objects $X_{\rho}$ and $Y_{\rho}$ of $\topos[1]$ such that there is canonical isomorphism:

\[ l^{2}(X_{\rho}) \simeq \rho(l^{2}(X)) \]
\[ l^{2}(Y_{\rho}) \simeq \rho(l^{2}(Y)) \]

Which are compatible with the canonical maps $\Delta_X:l^{2}(X) \rightarrow l^{2}(X) \otimes l^{2}(X)$ and $\Delta_Y:l^{2}(Y) \rightarrow l^{2}(Y) \otimes l^{2}(Y)$. Moreover $X_{\rho}$ is canonically defined, for example as the set of $v \in \rho(l^{2}(X))$ such that $\Vert v \Vert =1$ and $\Delta(v)= v \otimes v$.

The main result of this section is that:

\Prop{For any map $f:X\rightarrow Y$, there exists a unique map $f_{\rho}: X_{\rho} \rightarrow Y_{\rho}$ such that 

 \[ \Gamma_{f_{\rho}} = \rho(\Gamma_f) .\]
 
Moreover this defines a functor from the category of separating objects of $\topos$ to $\topos[1]$.}

The proof of this proposition will be concluded in \ref{proofCtrMorph}
}

\block{\label{functorialityLemma}Before proving the proposition we observe the following lemma:

\Lem{For any $f:X \rightarrow Y$ and $g:Y \rightarrow Z$ two maps between decidable objects of a topos, one has the following identities in the monoidal category of its Hilbert spaces:

\begin{equation} \label{functorialityLemmaEq1}
\Gamma_f ^{*} \Gamma_f = id_{l^{2}(X)}
\end{equation}
\begin{equation} \label{functorialityLemmaEq2} (\Delta_X \otimes id_{l^{2}(Y)} ) \circ \Gamma_f = (id_{l^{2}(X)} \otimes \Gamma_f) \circ \Delta_X
\end{equation}
\begin{equation} \label{functorialityLemmaEq3}
(id_{l^{2}(X)} \otimes \Delta_Y) \circ \Gamma_f = (\Gamma_f \otimes id_{l^{2}(Y)}) \circ \Gamma_f \end{equation}
\begin{equation} \label{functorialityLemmaEq4} \Gamma_{g \circ f} = ({\Gamma_f}^{*} \otimes id_{l^{2}(Z)}) \circ (id_{l^{2}(X)} \otimes \Gamma_g) \circ \Gamma_f 
\end{equation}

}

\Dem{Each of these relations is immediately checked internally in the topos on generators.}
}

\block{\label{functorialityProp}\Prop{Let $X$ and $Y$ be two decidable objects of a topos, and let $\Gamma : l^{2}(X) \rightarrow l^{2}(X) \otimes l^{2}(Y)$ be an operator which satisfies equations (\ref{functorialityLemmaEq1}), (\ref{functorialityLemmaEq2}) and (\ref{functorialityLemmaEq3}) of lemma \ref{functorialityLemma}. Then there exists a unique map $f:X \rightarrow Y$ such that $\Gamma= \Gamma_f$. }

\Dem{We work internally in the topos. As we will prove the existence and the uniqueness of the map $f$ it will indeed corresponds to a unique external map. The uniqueness is in fact automatic, if $\Gamma_f= \Gamma_g$ then for all $x$, $e_x \otimes e_{f(x)}= e_x \otimes e_{g(x)}$ which imply that $f(x)=g(x)$. We just have to prove the existence.

\bigskip

We define for $x,x' \in X$ and $y \in Y$:

\[ \gamma_{x' y}^{x} := \scal{\Gamma(e_x)}{ e_{x'}\otimes e_y}.\]

Let also $x,a,b \in X$ and $u,v \in Y$. 

On one hand:
\[ \scal{(\Delta_X \otimes id ) \circ \Gamma(e_x)}{e_a \otimes e_b \otimes e_u} = \delta_{a,b} \scal{\Gamma(e_x)}{e_a \otimes e_u} = \delta_{a,b} \gamma^{x}_{a,u} \]

on the other hand:
\[ \scal{(id \otimes \Gamma) \circ \Delta_X(e_x)}{e_a \otimes e_b \otimes e_u}= \scal{e_x \otimes e_x}{e_a \otimes \Gamma^{*}(e_b \otimes e_u)} =\delta_{x,a} \gamma^{x}_{b,u} \]

Hence, as $\Gamma$ satisfies relation (\ref{functorialityLemmaEq2}):

\[ \delta_{x,a} \gamma^{x}_{b,u} =  \delta_{a,b} \gamma^{x}_{a,u} \]

taking $x=a$ one obtains:

\[ \gamma^{a}_{b,u} = \delta_{a,b} \gamma^{a}_{a,u} \]

We will denote $\gamma^{a}_u$ for $\gamma^{a}_{a,u}$. One can deduce from this that:

\[ \scal{\Gamma(e_x)}{e_a \otimes e_v} = \gamma^{x}_{a, v} =\delta_{x,a} \gamma^{x}_v= \scal{e_x}{e_a} \gamma^{x}_v \]

As this holds for any $a \in X$ one obtains by linearity and continuity that for any $m \in l^{2}(X)$ one has:

\[  \scal{\Gamma(e_x)}{m \otimes e_v} = \scal{e_x}{m} \gamma^{x}_v \]

We will now use relation (\ref{functorialityLemmaEq3}). On one hand:

\[ \scal{(id \otimes \Delta_Y) \circ \Gamma (e_x)}{e_a \otimes e_u \otimes e_v} = \delta_{u,v} \scal{\Gamma(e_x)}{ e_a \otimes e_u}= \delta_{u,v} \delta_{x,a} \gamma^{x}_{u} \]

On the other hand:

 \begin{multline*}\scal{(\Gamma \otimes id) \circ \Gamma(e_x)}{e_a \otimes e_u \otimes e_v}  = \scal{\Gamma(e_x)}{\Gamma^{*}(e_a \otimes e_u) \otimes e_v} \\ = \scal{e_x}{\Gamma^{*}(e_a \otimes e_u)} \gamma^{x}_v=\delta_{x,a} \gamma^{x}_u \gamma^{x}_v \end{multline*}

Hence,

\[ \delta_{u,v} \delta_{x,a} \gamma^{x}_{u} = \delta_{x,a} \gamma^{x}_u \gamma^{x}_v \]

Hence (taking $x=a$ and $u=v$) one obtains that $(\gamma^{x}_{u})^{2} = \gamma^{x}_{u}$ hence for all $x,u$ one has $\gamma^{x}_u = 0$ or $\gamma^{x}_u = 1$, also, for all $ u \neq v$ one has $\gamma^{x}_u \gamma^{x}_v=0$ hence for each $x$ there is at most one $v$ such that $\gamma^{x}_v$ is non-zero and hence equal to $1$.

Finally as, by relation (\ref{functorialityLemmaEq1}), $\Gamma^{*}\Gamma =1$, the norm of $\Gamma(e_x)$ is one, there is at least one $u$ such that $\gamma^{x}_u=1$. Hence for each $x \in X$ there exists a unique $u$ such that $\gamma^{x}_u=1$, hence there is a map $f$ defined by $f(x)=u$ and because:

\[ \scal{\Gamma(e_x)}{ e_{x'} \otimes e_y} = \gamma^{x}_{x',y} = \delta_{x,x'} \gamma_y^{x} = \delta_{x,x'} \delta_{f(x),y} \]

one indeed has $\Gamma(e_x)=e_x \otimes e_{f(x)}$ i.e. $\Gamma = \Gamma_{f}$.

}

}

\block{\label{proofCtrMorph}We can now prove proposition \ref{functorial}:

\Dem{As $\rho$ is monoidal, $\rho(\Gamma_f)$ satisfies all the relations of the lemma \ref{functorialityLemma}, hence by proposition \ref{functorialityProp} one obtain a unique map $f_{\rho}$ such that $\rho(\Gamma_f) = \Gamma_{f_{\rho}}$. The functoriality comes immediately from relation (\ref{functorialityLemmaEq4}) of lemma \ref{functorialityLemma} and the fact that as $\rho$ is monoidal it preserves this relation.
}

}

\subsection{Construction of the Geometric morphism}
\label{secCtrGeomMorph}

\blockn{In the previous subsection one has obtain from the representation $\rho$ a functor $F_{\rho}:X \mapsto X_{\rho}$ from the category of separating objects of $\topos$ to the the category of all objects of $\topos[1]$ which in some sense is compatible with $\rho$. As we assume that $\topos$ is locally separated, the category of separating object of $\topos$ endowed with the restriction of the canonical topology of $\topos$ is a site of definition for $\topos$. Moreover, this site has all finite limits (except the empty limit) and arbitrary (small) co-products, hence in order to check that the functor $F_{\rho}:X \mapsto X_{\rho}$ we have constructed extend into the $f^{*}$ part of a geometric morphism from $\topos[1]$ to $\topos$ we essentially\footnote{see the proof of \ref{CtrGeomMorph}.} need to check that it preserves finite limits, arbitrary co-products and send epimorphisms in $\topos$ to epimorphisms in $\topos[1]$.}

\block{\label{Fsurjection1}\Lem{Let $X$ and $Y$ be two separating objects of $\topos$ and $f:X\twoheadrightarrow Y$ an epimorphism such that $X$ has a bounded cardinal in $\Tcal_{/Y}$. Then $F_{\rho}(f)$ is an epimorphism in $\topos[1]$.}

\Dem{As $f$ has finite fibres of bounded cardinal there exists an operator $g:l^{2}(X) \rightarrow l^{2}(Y)$ such that $g(e_x)=e_{f(x)}$ (the norm of this map is smaller than the square root of the bound of the size of the fibres). In particular,

\[ \Gamma_f = (Id \otimes g) \circ \Delta_X \]

Moreover, as $f$ is a surjection, if one defines:

\[ h(e_y) = \frac{1}{n} \sum_{f(x)=y} e_x \]

where $n$ denotes the number of pre-images of $y$ by $f$, then $g \circ h = Id$.

Interpreting this in $\topos[1]$, the identity $\Gamma_f = (Id \otimes g) \circ \Delta_X $ characterize $\rho(g)$ as the map which send $e_x$ to $e_{f_{\rho}(x)}$ and the existence of the map $\rho(h)$ such that $\rho(g) \rho(h) = Id$ show that $f_{\rho}$ is surjective.
}

}

\block{\label{Fsurjection2}\Lem{Let $f:X \twoheadrightarrow Y$ be an arbitrary surjection between two separating objects of $\topos$ then $f_{\rho}$ is also a surjection.}

\Dem{By \ref{SepImpFinitness} applied to $\Tcal_{/Y}$, the object $X$ admit a covering by a net of subobjects $U \subset X$ such that the restriction of $f$ to $U$ has bounded cardinal in $\Tcal_{/Y}$. For each $U$ of this net, let $Y_U$ be the image of $U$ by $f$. From the previous lemma, the image by $F_{\rho}$ of the maps from $U$ to $Y_U$ are all surjective, hence all we have to prove is that the family of $F_{\rho}(Y_U)$ cover $F_{\rho}(Y)$. Now, as in the proof of lemma \ref{Fsurjection1} for each $U \subset X$ the corresponding map $i_U : l^{2}(Y_U) \rightarrow l^{2}(Y)$ is send by $\rho$ on the operator corresponding to the map from $F_{\rho}(Y_U)$ to $F_{\rho}(Y)$. The fact that the $Y_U$ cover $Y$ is translated into the fact that $i_U i_U^{*}$ converge weakly to the identity and hence is transported by $\rho$, which concludes the proof.}
}

\block{\Lem{$F_{\rho}$ preserve arbitrary\footnote{Only those indexed by decidable sets in fact, but we are not working internally in $\Ecal$ any more and we are assuming the law of excluded middle in the base topos, so this does not really matter.} co-products.}

\Dem{One has $U=\coprod_i U_i$, if and only if the $U_i$ are complemented sub-objects of $U$ with $t_i : l^{2}(U_i) \rightarrow l^{2}(U)$ the corresponding map, such that:

\[t_i^{*} t_i =Id_{l^{2}(U_i)} \]
\[t_i^{*} t_j =0 \text{ when $i \neq j$.} \]
\[\sum_i t_i t_i^{*} = Id_{l^{2}(U)} \]
where the infinite sum in the last equality mean the supremum of the the net of finite sum.

All those properties are clearly preserved by $\rho$, hence the coproduct is preserved by $F_\rho$.
}
}

\block{\Lem{$F_\rho$ preserves binary products and equalizers of pair of maps.}

\Dem{ One has $l^{2}(X \times Y) \simeq l^{2}(X) \otimes l^{2}(Y)$ and through this isomorphism, the operator $\Delta$ for $X \times Y$ is the tensor product of those for $X$ and $Y$, and this completely characterize the space $l^{2}(X \times Y)$ and its operator $\Delta$. As $\rho$ is monoidal it preserve this ``characterization" of the product and hence $F_\rho$ preserve binary product.

\bigskip 

Let
\[X_0 \overset{i}{\hookrightarrow} X \overset{f,g}{\rightrightarrows} Y\]

be the equalizer of a pair of maps between separating objects of $\topos$. Then for all $x\in X$ one has:

\[\begin{array}{r c l l} \Gamma_f ^{*} \Gamma_g (e_x)  &= &e_x &\text{ if $f(x)=g(x)$ } \\ & =& 0 & \text{ otherwise} \end{array} \]

Hence $\Gamma_f ^{*} \Gamma_g$ is a projection and $l^{2}(X_0)$ is isomorphic to $\Gamma_f ^{*} \Gamma_g (l^{2}(X))$. This is also preserved by $\rho$, hence $F_{\rho}$ preserves equalizer.
}
}

\block{\label{CtrGeomMorph}\Prop{There exists a geometric morphism $\chi_{\rho} : \topos[1] \rightarrow \topos$ such that $F_{\rho}$ is the restriction of $\chi_{\rho}^{*}$ to the subcategory of separating objects. }

\Dem{As the subcategory of separating object form a generating familly of $\topos$ it is a site of definition when we equiped it with the canonical topology of $\topos$. Hence, all we have to prove is that $F_{\rho}$ is flat and continuous for the canonical topology of $\topos$. It is continuous because it preserves surjections and arbitrary co-products. We just have to prove that it is flat i.e. that (internally in $\topos[1]$) :

\begin{itemize}
\item $\exists x \in F_{\rho}(X)$ for some objects $X$.

\item If $c_1 \in F_{\rho}(X_1)$ and $c_2 \in F_{\rho}(X_2)$ then there exists an object $X_3$ with two arrow $f_i : X_3 \rightarrow X_i$ and a $c_3 \in F_{\rho}(X_3)$ such that $f_i x_3 = x_i$.

\item If $c \in F_{\rho}(X)$ and there is two arrow $f,g : X \rightrightarrows Y$, such that $f.c = g.c$ then there is an arrow $h:X' \rightarrow X$ and a $c' \in F_{\rho}(X')$ such that $h.c' = c$ and $f.h=g.h$.

\end{itemize}

The second condition follows from the fact that $F_{\rho}$ commute to binary product, and the third from the fact that it commute to equalizer. The first condition is generally a consequence of the fact that the functor preserve the terminal object, but in our situation if $\topos$ is not separated there is no terminal object in our site. Instead, $F_{\rho}(X_0)$ is inhabited because $\rho(l^{2}(X_0))$ contains a vector of strictly positive norm, which concludes the proof.
}

}

\block{Finally, if $\mu : \rho \rightarrow \rho'$ is a morphism in the sense of definition \ref{Defmorphofrep}, then $\mu$ induce a natural transformation between $F_{\rho}$ and $F_{\rho'}$ because for any separating object $X$ of $\Tcal$, an $x \in \rho(l^{2}(X))$ such that $\Delta(x) = x \otimes x$ and $\Vert x \Vert =1$ will be send by $\mu$ on an $x$ with the same properties in $\rho'(l^{2}(X))$ hence inducing a map from $F_{\rho}(X)$ to $F_{\rho'}(X)$ and this is natural in $X$. And a natural transformation between the $F_{\rho}$ extend uniquely into a natural transformation between the induced geometric morphism

hence we have constructed a functor from (non-degenerate) representation of $\Hcal(\Tcal)$ in $\Ecal$ to geometric morphisms from $\Ecal$ to $\Tcal$.

}

\subsection{Proof of the ``reduced" theorem \ref{mainThRed} }
\label{ProofmainThred}

\block{If we start from a geometric morphism $\chi: \topos[1] \rightarrow \topos$ then we associate to it a representation $\rho=\chi^{\sharp}$, defined on all of $\Hcal(\topos)$ and we can restrict it to $\Hcal^{red}(\topos)$.

As for any object $X$ of $\topos$, one has $\chi^{\sharp}(l^{2}(X)) \simeq l^{2}(\chi^{*}(X))$, and for any map $f$ one has $\chi^{\sharp}(\Gamma_f)= \Gamma_{\chi^{*} f}$ (up to te canonical isomorphism), it is clear that the geometric morphism reconstructed from this representation $\rho$ will be isomorphic to $\chi$ and that there will be a bijection between morphisms of representations and natural transformation of geometric morphisms. So all we have to do to finish the proof of the theorem is the to prove the following: }

\block{\label{proofcompatibilityRED}\Prop{For any reprensation $\rho$ of $\Hcal^{red}(\topos)$ in a topos $\topos[1]$, $\rho$ is equivalent to the representation $\chi_{\rho}^{\sharp}$ induced by the geometric morphism $\chi_{\rho}:\topos[1] \rightarrow \topos$ constructed out of $\rho$. }

\newcommand{\rhop}{\chi_{\rho}^{\sharp}}

\Dem{We will first prove that $\rho$ and $\rhop$ are equivalent on the full-subcategory of Hilbert spaces of the form $l^{2}(X)$ for $X$ a separating object of $\topos$. The geometric morphism $\chi_{\rho}$ has been constructed to full-fill the relation $l^{2}(\chi^{*}(X)) \simeq \rho(l^{2}(X))$ hence one has indeed an isomorphism $\rho(l^{2}(X)) \simeq \rhop(l^{2}(X))$ for any separating object $X$ of $\topos$. Also, by construction, these isomorphisms are compatible with the maps $\Delta_X$ and  $\Gamma_f$ for any separating object $X$ of $\topos$ and any maps $f:Y \rightarrow X$.

\bigskip

Moreover, if $S \subset X$, then $\rho$ and $\rhop$ also agree on the projection $P$ on $l^{2}(X)$ corresponding to the multiplication by the characteristic function of $S$. As any operator of multiplication by a bounded function on $X$ is a supremum of linear combination of such projection, $\rho$ and $\rhop$ also agree on operators on $l^{2}(X)$ defined by multiplication by complex bounded functions on $X$.

\bigskip

Let now $f :l^{2}(X) \rightarrow l^{2}(Y)$ be any bounded operator in any topos. Let $m$ the function of matrix element of $f$ on $X \times Y$, ie internally, $m(x,y):=\scal{e_y}{f(e_x)}$. The operator $M_m$ of multiplication by $m$ on $l^{2}(X) \otimes l^{2}(Y)$ is characterized by the relation :

\[ M_m = ( id_{l^{2}(X)} \otimes \Delta_Y^{*}) \circ (id_{l^{2}(X)} \otimes f \otimes id_{l^{2}(Y)}) \circ (\Delta_X \otimes id_{l^{2}(Y)})  \]

which is easily checked internally on generators.

In particular, if $X$ and $Y$ are two separating objects of $\topos$ and $f :l^{2}(X) \rightarrow l^{2}(Y)$ is any bounded operator, then the operator of multiplication by matrix elements of $\rho(f)$ is $\rho(M_m)$ because of the previous relation is preserved by $\rho$. Hence as $\rho(M_m)= \rhop(M_m)$ one can deduce that $\rho(f)$ and $\rhop(f)$ have the same matrix elements and hence are equals.

\bigskip

As $l^{2}(X)$ for $X$ a separating bound is a generator of $\Hcal^{red}(\Tcal)$,this is enough to prove that the two $*$-functors $\rho$ and $\rhop$ are isomorphic on the whole category $\Hcal^{red}(\Tcal)$ by the result mentioned in section $2$.
}

 }

\section{On the category $\Hcal(\Tcal)$ and its representations}
\label{secProofUnred}
\subsection{The category $\Hcal(\Tcal_{/X})$}

\block{In this subsection, $\Tcal$ is an arbitrary topos, and $X$ is a decidable object of $\topos$. As previously, this defines an object $l^{2}(X)$ of $\Hcal(\Tcal)$ which is endowed with a multiplication $\Delta^{*} : l^{2}(X) \otimes l^{2}(X) \rightarrow l^{2}(X) $.}

\block{If $\Hcal$ is a Hilbert space of $\Tcal_{/X}$, then, internally in $\Tcal$ it corresponds to a $X$-indexed family of Hilbert spaces $(\Hcal_x)_{x \in X}$. As $X$ is decidable one can construct the orthogonal sum of such a family:

\[ \Sigma \Hcal := \oplus_{x \in X} \Hcal_{x} \] 

which is a Hilbert space of $\Tcal$, and which comes endowed with an ``action" of $l^{2}(X)$:

\[ \begin{array}{c c c c}
M:  &l^{2}(X) \otimes \Sigma \Hcal &\rightarrow & \Sigma\Hcal \\
& e_y \otimes (h_x)_{x \in X} & \mapsto & (\delta_{x,y} h_x)_{x \in X} 
\end{array}
\]

which admit an adjoint $M^{*}((h_x)_{x \in X}) = \sum_{a}\left( e_a \otimes (\delta_x^{a} h_x)_{x \in X} \right)$, hence is indeed a morphism in $\Hcal(\Tcal)$ and satisfies $MM^{*}=Id_{\Sigma \Hcal}$ (and it is indeed an action because it is the multiplication of a sequence of vector by a sequence of scalar).
 
Finally, these structures satisfy the following additional relation:
\[M^{*}M = (id_{\Sigma \Hcal} \otimes \Delta^{*}) \circ (M^{*}\otimes id_{l^{2}(X)}) \]

i.e. the following diagram commute:
\[
\begin{tikzcd}[ampersand replacement=\&,scale=2]
\Sigma \Hcal \otimes l^{2}(X) \arrow{rr}{M^{*} \otimes id_{l^{2}(X)}} \arrow{d}{M} \& \& \Sigma \Hcal \otimes l^{2}(X) \otimes l^{2}(X) \arrow{d}{id_{\Sigma \Hcal} \otimes \Delta^{*}} \\
\Sigma \Hcal \arrow{rr}{M^{*}} \&\&  \Sigma \Hcal \otimes l^{2}(X)\\
\end{tikzcd}
\]

Or to put it another way, $M^{*}$ is $l^{2}(X)$-linear when $\Sigma \Hcal \otimes l^{2}(X)$ is endowed with the action of $l^{2}(X)$ on the second component of the tensor product.

This is also easily checked internally on an element of the form $(h_x) \otimes e_y$.

}

\block{\label{Defl2module}\Def{A $l^{2}(X)$-module is an object $H \in \Hcal(\Tcal)$ endowed with an action $M: H\otimes l^{2}(X) \rightarrow H$ of the monoide $(l^{2}(X),\Delta^{*})$, which satisfies the two additional relation: $MM^{*}=id_H$ and $M^{*}M = (id_{H} \otimes \Delta^{*}) \circ (M^{*}\otimes id_{l^{2}(X)})$.

A morphism of $l^{2}(X)$-module is a morphism of module in the usual sense. 
}

We already mentioned our main example: $\Sigma \Hcal$ is a $l^{2}(X)$-module, and one can check that if $(f_x)$ is a globally bounded family of operators $f_x:\Hcal_x \rightarrow \Hcal'_x$, then $\Sigma f$ is a morphism of $l^{2}(X)$-modules from $\Sigma \Hcal$ to $\Sigma \Hcal'$.}

\block{\label{l2Xmodeqtoslice}
\Prop{The following three categories are equivalent:

\begin{itemize}

\item The category of $l^{2}(X)$-modules, as in definition \ref{Defl2module}.

\item The category of non-degenerate representations of $C_0(X)$. I.e. the category of Hilbert spaces $H \in \Hcal(\Tcal)$ endowed with an action of $C_0(X)$ making them, internally in $\Tcal$, a non-degenerate\footnote{Non-degenerate meaning here that the image of the action map $C_0(X) \times H \rightarrow H$ spam a dense subspace of $H$.} $*$-representation of $C_0(X)$.

\item $\Hcal(\Tcal_{/X})$

\end{itemize}
}

\Dem{The equivalence between the second and the third category is proved internally and is extremely classical (and is constructive if we assume $X$ decidable): any non-degenerate representation $H$ of $C_0(X)$ decompose into a direct sum of $H_x = e_x(H)$ where $e_x$ is the characteristic function of $\{x\}$.

Suppose we start with a $l^{2}(X)$-module $H$. It is in particular a Hilbert space endowed with an action of the algebra $(l^{2}(X),*)$. As we mentioned, the axiom $M^{*}M = (id_{H} \otimes \Delta^{*}) \circ (M^{*}\otimes id_{l^{2}(X)})$ mean that the map $M^{*}: H \rightarrow H \otimes l^{2}(X) $ is $l^{2}(X)$-linear when we endow $H \otimes l^{2}(X)$ with the action of $l^{2}(X)$ on the second component only, and hence, $M$ and $M^{*}$ exhibit $H$ as a $l^{2}(X)$-linear retract of $H \otimes l^{2}(X)$. But the action of $l^{2}(X)$ on $H \otimes l^{2}(X)$ clearly extend into a non-degenerate representation of $C_0(X)$ and hence (as $M^{*}M$ is also $C_0(X)$-linear) the action of $l^{2}(X)$ on $H$ also extend into a non-degenerate representation of $C_0(X)$.

Conversely, by the equivalence of the second and the third categories, any $C_0(X)$ module is up to unique isomorphism of the form $\Sigma \Hcal$ and hence comes from a $l^{2}(X)$-module. This already proves the equivalence at the level of objects, but, as $l^{2}(X)$ is internally dense in $C_0(X)$ for the uniform norm, a map is $C_0(X)$-linear if and only if it is $l^{2}(X)$-linear and hence this concludes the proof.
}

}

\block{\label{pullbackl2module}We conclude this subsection by investigating how, if $f :X \rightarrow Y$ is any map in $\Tcal$, the ``pullback along $f$" functor from $\Hcal(\Tcal_{/Y})$ to $\Hcal(\Tcal_{/X})$ can be expressed in terms of $l^{2}(Y)$ and $l^{2}(X)$-modules and the operator $\Gamma_f : l^{2}(X) \rightarrow l^{2}(X) \otimes l^{2}(Y)$ defined in \ref{DefGammaf}.

Let $V$ be an $l^{2}(Y)$-module, consider the following map:

\[ P_{V,f} : V \otimes l^{2}(X) \overset{\Gamma_f}{\rightarrow} V \otimes l^{2}(X) \otimes l^{2}(Y) \overset{M_V}{\rightarrow} V \otimes l^{2}(X) \]

where the first map is $\Gamma_f$ acting on the second component and the second map is the structural multiplication of $l^{2}(Y)$ on $V$ after re-ordering of the terms.

\Prop{\begin{itemize}

\item $P_{V,f}$ is a $l^{2}(X)$-linear orthogonal projector of $V \otimes l^{2}(X)$. In particular, its image $P_{V,f}(V \otimes l^{2}(X) )$ is a $l^{2}(X)$-module. 

\item If, as a $l^{2}(Y)$-module, one has $\displaystyle V = \sum_{y \in Y} V_y$, then:

\[ P_{V,f}(V \otimes l^{2}(X)) \simeq \sum_{x \in X} V_{f(x)}\]

as a $l^{2}(X)$-module.

\item If $h:V \rightarrow V'$ is a map of $l^{2}(Y)$-module, then $(h \otimes Id_{l^{2}(X)})$ commutes to the $P_{V,f}$ and induce the pullback of $h$ along $f$ between $P_{V,f}(V \otimes l^{2}(X))$ and $P_{V',f}(V' \otimes l^{2}(X) )$.
\end{itemize}
}

\Dem{Everything can be checked internally on generators, indeed:

\[ P_{V,f} ( v \otimes e_x) = v.e_{f(x)} \otimes e_x \]

where the point denotes the action of $l^{2}(Y)$ on $V$, hence one immediately obtains that $P_{V,f} P_{V,f} (v \otimes e_x ) = v.e_{f(x)} \otimes e_x = P_{V,f}(v \otimes e_x)$ and that for any $x' \in X$:

\[ P_{V,f}(v \otimes e_x) .e_{x'} = v.e_{f(x)} \otimes e_x.e_x' = P_{V,f}(v \otimes e_x) \delta_{x,x'} = P_{V,f}(v \otimes e_x.e_x'). \] 

This concludes the first point, for the second, if $V = \sum_{y \in Y} V_y$ and if $v \in V_y$, then $P_{V,f}(v\otimes e_x)$ is $(v \otimes e_x) \delta_{f(x),y}$, hence $P_{V,f}(V \otimes l^{2}(X)$ is the subspace generated by the $(v \otimes e_x)$ for $v \in V_{f(x)}$ which is isomorphic to $\sum_{x \in X} V_{f(x)}$ (they have the same generators with the same scalar product between them).

For the third point, $h(v \otimes e_x) = h(v) \otimes e_x$ and $h$ is $l^{2}(Y)$ linear, hence:

\[  P_{V',f}(h(v) \otimes e_x) = h(v).e_{f(x)} \otimes e_x = h(v.e_{f(x)}) \otimes e_x = h(P_{V,f}(v \otimes e_x))  \]

And the action of $h$ on $P_{V,f}(v \otimes e_x)$ send $v \otimes e_x$ for $v \in V_{f(x)}$ to $h(v) \otimes e_x$ with $h(v) \in V'_{f(x)}$ hence it is indeed the pullback of $h$ along $f$.

}

}

\subsection{Tensorisation by square integrable Hilbert space}

\block{\label{LemmaSepimpSI}\Lem{Let $\Tcal$ be a separated boolean topos, then every Hilbert spaces of $\Tcal$ is square integrable.}

\Dem{Let $\Hcal$ be a Hilbert space of $\Tcal$. As $\Tcal$ is separated, theorem \ref{SepImpFinitness} implies it is generated by objects of (finite) bounded cardinal hence $\Hcal$ can be covered by maps of the form $S \rightarrow \Hcal$ with $S$ of bounded cardinal. Such a map from $S$ to $\Hcal$ can be extended by linearity into a globally bounded map from $l^{2}(S)$ to $\Hcal$ (whose norm will be smaller than $\sqrt{n}$ for any $n$ such that $|S| \leqslant n $). Hence $\Hcal$ is covered by maps from square integrable Hilbert spaces and hence is itself square integrable.}
}

\block{\label{tensorbySI}\Prop{Let $\Tcal$ be a boolean locally separated topos, $\Hcal$ a $\Tcal$-Hilbert space and $\Hcal_i$ a square integrable $\topos$-Hilbert space. Then $\Hcal \otimes \Hcal_i$ is square integrable.}

\Dem{$\Hcal_i$ being square integrable it can be covered by map from $l^{2}(X)$ with $X$ separating, hence $\Hcal \otimes \Hcal_i$ can be covered by maps of the form $\Hcal \otimes l^{2}(X)$ for $X$ separating. Hence it is enough to prove the result in the case $\Hcal_i = l^{2}(X)$ for $X$ a separating object.

Let $p:\Tcal_{/X} \rightarrow \Tcal$ the canonical map. Internally in $\Tcal$, a $\Tcal_{/X}$ Hilbert space is just a $X$-indexed family of Hilbert spaces, $p^{*}(\Hcal)$ is the constant family equal to $\Hcal$, and $l^{2}(X) \otimes \Hcal$ is the orthogonal sum of this family. 

$\Tcal_{/X}$ being separated, $p^{*}(\Hcal)$ is a square integrable Hilbert space in $\Tcal_{/X}$ by lemma \ref{LemmaSepimpSI}, hence it admit a covering by globally bounded maps $l^{2}(Y_i) \rightarrow p^{*}(\Hcal)$ for $Y_i \rightarrow X$ objects of $\Tcal_{/X}$. Taking the orthogonal sum along $X$ internally in $\Tcal$ one obtains a series of maps $l^{2}(Y_i)\rightarrow \Hcal \otimes l^{2}(X)$ which also form a covering, and hence $\Hcal \otimes l^{2}(X)$ is square integrable.
}

}

\subsection{Proof of the ``unreduced" theorem \ref{mainThunred}.}

\block{Let first $\rho$ be a representation of $\Hcal(\Tcal)$ in $\Ecal$ in the sense of definition \ref{DefRepUnred}. Let $\rho^{red}$ be the restriction of $\rho$ to the full subcategory $\Hcal^{red}(\Tcal) \subset \Hcal(\Tcal)$.  Clearly, $\rho^{red}$ is a representation of $\Hcal^{red}(\Tcal)$ in the sense of definition \ref{DefRepReduced}.
In particular, by theorem \ref{mainThRed}, which has been proved in \ref{ProofmainThred}, there exists a geometric morphism $f:\Ecal \rightarrow \Tcal$ such that $\rho^{red}$ is isomorphic to $f^{\sharp}$ (as symmetric monoidal $C^{*}$-functor) on $\Hcal^{red}(\Tcal)$. 

\Prop{The isomorphism between $\rho^{red}$ and $f^{\sharp}$ on $\Hcal^{red}(\Tcal)$ extend canonically\footnote{In fact, uniquely as the characterization of the isomorphism given during the proof will show.} into an isomorphism of symmetric monoidal $C^{*}$-functor between $\rho$ and $f^{\sharp}$ on $\Hcal^{red}(\Tcal)$.}

\Dem{Let $H$ be any object of $\Hcal(\Tcal)$, let $X$ be any inhabited separating object of $\Tcal$ and let $H_0 = l^{2}(X) \in \Hcal^{red}(\Tcal)$. One already knows that $l^{2}(f^{*}(X)) \simeq f^{\sharp}(H_0) \simeq \rho(H_0)$. By proposition \ref{tensorbySI}, the $\Tcal$-Hilbert space $H \otimes H_0$ is in $\Hcal^{red}(\Tcal)$, hence one has an isomorphism $\rho(H \otimes H_0) \simeq f^{\sharp}(H \otimes H_0)$, and even an isomorphism $\rho(H) \otimes l^{2}(f^{*}X) \simeq f^{\sharp}(H) \otimes l^{2}(f^{*}X))$, moreover, as the isomorphism between $\rho$ and $f^{\sharp}$ on $\Hcal^{red}$ is functorial and symmetric monoidal, this last isomorphism is an isomorphism of $l^{2}(f^{*}X)$-module in the sense of definition \ref{Defl2module}. In particular, internally in $\Ecal$ it comes from a family $(\lambda_x)_{x \in f^{*}X}$ of isomorphisms between $\rho(H)$ and $f^{\sharp}(H)$.

If now $X$ and $Y$ are two inhabited separating objects of $\Tcal$ and $g:X\rightarrow Y$ any map. The previous construction gives us (internally in $\Ecal$) two families of isomorphisms $(\lambda_x)_{x \in f^{*}X}$ and $(\mu_y)_{y \in f^{*}Y}$ between $\rho(H)$ and $f^{\sharp}(H)$.

The description given in \ref{pullbackl2module} of pullbacks of $l^{2}(Y)$-modules (and their morphisms) into $l^{2}(X)$-modules (and their morphisms) along $g$ is preserved by $\rho$ and $f^{\sharp}$ just because they are monoidal symmetric, hence one obtains in $\Hcal(\Tcal)$ exactly the relation asserting the that the family of isomorphisms $(\lambda_x)$ is the pullback of the family $(\mu_y)$, i.e. for all $x\in X$, $\lambda_x = \mu_f(x)$.

Applying this observation to the two projections $X \times X \rightarrow X$ one obtain that for all $x,x' \in X \times X$ , $\lambda_x = \lambda_{x,x'} = \lambda_x'$. Hence the isomorphism $\lambda_x$ does not depend on $x$ and as $X$ is inhabited this proves (in $\Ecal$) that there is a canonical isomorphism between $\rho(H)$ and $f^{\sharp}(H)$ and hence this gives an external isomorphism $\mu_H : \rho(H) \rightarrow f^{\sharp}(H)$ (which does not depend on $X$ neither). This isomorphism is entirely characterized by the fact that for some (or any) inhabited separating object $X$ of $\Tcal$ the isomorphism $f^{\sharp}(H) \otimes l^{2}(f^{*}X) \simeq f^{\sharp}(H \otimes l^{2}(X)) \simeq \rho(H \otimes l^{2}(X)) \simeq \rho(H) \otimes l^{2}(f^{*}X) $ is equal to $\mu_H \otimes Id_{l^{2}(f^{*}X)}$.

The fact that this $\mu_H$ is indeed a monoidal symmetric natural transformation then follow simply from this characterization, and the fact that the isomorphism between $f^{\sharp}$ and $\rho$ on square integrable Hilbert spaces is itself functorial and symmetric monoidal. 
}
}

\block{At this point we have proved that the two constructions forming the equivalence in theorem \ref{mainThunred} are inverse of each other at the level of object (up to canonical isomorphism).  To completely conclude the proof of this theorem one need to check that the notion of morphisms on the two side are the same.

let $\phi: \rho \rightarrow \rho'$ be an isometric inclusion of symmetric monoidal representation (in the sense of definition \ref{Defmorphofrep}  ), then $\phi$ restricted to $\Hcal^{red}(\Tcal)$ is a morphism between $\rho^{red}$ and $\rho'^{red}$ and hence induce a natural transformation between the corresponding geometric morphism $f$ and $f'$ that we will denote by $F(\phi)$. Conversely, any natural transformation $\mu$ between $f$ and $f'$ will induce an inclusion of $f^{\sharp} \simeq \rho$ into $f'^{\sharp} \simeq \rho'$ which we will denote by $G(\mu)$. It is immediate that $F (G(\mu)) = \mu$ (up to the canonical isomorphism). And it is already proved that the two morphisms of representations $G(F(\phi))$ and $\phi$ agree when restricted to $\Hcal^{red}(\Tcal)$. The following proposition allows to conclude.

\Prop{If $\phi$ and $\phi'$ are two morphisms of representations from $\rho$ to $\rho'$ (two non-degenerate representations of $\Hcal(\Tcal)$ in $\Ecal$ ) which agree on objects of $\Hcal^{red}(\Tcal)$ then $\rho = \rho'$.}

\Dem{Let $H$ be any object of $\Hcal(\Tcal)$ and $H'$ be any object of $\Hcal^{red}(\Tcal)$ whose image by $\rho$ in $\Hcal(\Ecal)$ contains a vector of norm $1$.

Then as $\phi$ and $\phi'$ are monoidal one has $\phi_H \otimes \phi_{H'} = \phi_{H \otimes H'}$ and similarly for $\phi'$, but as $H'$ and $H \otimes H'$ are square integrable (by assumption for the first one and by proposition \ref{tensorbySI} for the second) $\phi$ and $\phi'$ are the same on them and hence $\phi_H \otimes \phi_{H'} = \phi'_H \otimes \phi_{H'}$. Then, internally in $\Ecal$, one can take some $h \in \rho(H')$ of norm one, and composing by the corresponding maps from $\C$ to $\rho(H')$ and from $\rho'(H')$ to $\C$ one obtains that internally $\phi_H = \phi'_H$ which immediately show that the same equality hold externally and concludes the proof of the proposition, and of theorem \ref{mainThunred}.

}
 
}

\section{Toward a generalized Gelfand duality ?}
\label{secToward}

\blockn{The goal of this last section is to informally discus the possibility to extend the results of this paper into some sort of ``non-commutative Gelfand duality" that would relates geometric (or topological) objects with objects from the world of operators algebras.}

\blockn{A pre-requisite to understand what follows, is to understand how toposes can be fully faithfully embedded into a category of localic stacks (in fact ``geometric localic stacks", hence localic groupoids) as clearly explained in \cite{bunge1990descent}. Roughly: To any topos $\Tcal$ one can associate the stacks defined by $\hat{\Tcal}(\Lcal) = \hom(\Lcal,\Tcal)$ where $\Lcal$ is any locale\footnote{In what follows we identify a locale with its topos of sheaves, or equivalently, we actually mean ``localic topos" when we write ``locale".} and where non-invertible natural transformations in $\hom(\Lcal,\Tcal)$ have been dropped. This defines a fully faithful embeddings of the $2$-category of toposes into the $2$-category of geometric\footnote{that is the stacks that can be represented by a groupoid.} stacks on the category of locales endowed with the open surjection topology. This embeddings is the right adjoint of the functor which send any localic groupoid to the category of equivariant sheaves over it. }

\blockn{We will first explain how the two construction involved in the present ``duality" can be defined at a very high level of generality, for example the category $\Hcal(\Tcal)$ is defined for any topos $\Tcal$ without any assumption on it, but we can do even better.}

\blockn{If we start from an arbitrary monoidal symmetric $C^{*}$-category $C$ one can try to define its ``spectrum" $\spec C$ as the ``classifying topos" of symmetric monoidal $*$-representations of $C$. This $\spec C$ has in general no reason to exist as a topos, but it always makes sense as a prestack over the category of toposes or of locales: }

\blockn{\Def{Let $C$ be a symmetric monoidal $C^{*}$-category, then $\spec C$ is the pre-stack over the category of locales defined by: $(\spec C) ( \Lcal)$ is the groupoid of symmetric monoidal representations of $C$ in $\Lcal$ and symmetric monoidal $*$-isomorphisms between them, and for any morphism of locale $f:\Lcal \rightarrow \Lcal'$, $(\spec C)(f)$ is just $f^{\sharp}$. 
}
}

\blockn{Defined this way, $\spec C$ is not a stack for a reasonable topology one the category of locales (except the topology of etale surjection). The problem is that Hilbert spaces, because they are not pulled-back as sheaves, does not satisfy reasonable descent properties. Fortunately, we have presented the solution to this problem in \cite{henry2014localic} :  we proved that (because Hilbert space are actually pulled-back as locales) localic Hilbert spaces do descend along open surjections (in fact along any morphism which is of effective descent for locales), and form in fact the stackification of the pre-stack of Hilbert spaces.

Hence we have to change the definition of $\spec C$ by replacing $\Hcal(\Lcal)$ by the the category of localic Hilbert space of $\Lcal$ as defined in \cite{henry2014localic}, and hence $\spec C$ classifies representations of $C$ on localic Hilbert spaces, and with this definition $\spec C$ will be a stack for the topology of open surjection (as well as for any topology composed of effective descent morphisms).

}

\blockn{Still thanks to results of \cite{henry2014localic}, one can define a stacks ``$\Hcal$" of localic Hilbert spaces. One can even prove that it is a geometric stacks : its diagonal is a localic map because of \cite[section 3.5]{henry2014localic} and it admit an open surjection from the locale classifying the pre-Hilbert structure on an infinite dimensional $\Q$-vector space.

Moreover, for any localic stacks $\Gcal$ one can define a symmetric monoidal $C^{*}$-category $\Hcal(\Gcal)$ whose objects are morphisms of stacks from $\Gcal$ to $\Hcal$ (hence, essentially continuous fields of Hilbert spaces over $\Gcal$).
}

\blockn{At this point, and up to some set theoretical issues that we do not want to talk about here, one obtains two functors ``$\spec$" and ``$\Hcal$" which essentially form an adjoint pair between localic stacks (for example for the topology of open surjections) and symmetric monoidal $C^{*}$-category.}

\blockn{That this adjoint pair that we think can provide this generalized duality. It is not reasonable to think that every localic stacks can be reconstructed out of its category of Hilbert spaces: For example, let $\set[\Ob^{+}]$ be the classifying topos for the theory of inhabited objects (i.e. the topos of preshaves over the opposite of the category of finite non-empty sets) then one can see that $\Hcal(\set[\Ob^{+}])$ is equivalent to the monoidal $C^{*}$-category of Hilbert spaces.

It is not reasonable neither to think that any localic stacks can actually be reconstructed exactly as $\Gcal \simeq \spec \Hcal \Gcal$. As we have seen in the present paper, one needs to restrict the notion of representations considered in the definition of $\spec C$ (in our case we only considered ``non-degenerate normal" representations). }

\blockn{But for example, one can easily prove using the same kind of techniques as in section \ref{secFunctoOnObj} and \ref{secCtrGeomMorph} that for any locally decidable topos $\Tcal$ (that is a topos such that every object can be covered by decidable objects), the stack associated to $\Tcal$ is a sub-stacks of $\spec \Hcal(\Tcal)$. Hence for any locally decidable topos $\Tcal$ there exists a notion of ``good" representation of $\Hcal(\Tcal)$ such that $\Tcal$ is the classifying topos for these ``good" representations.

It opens several questions (which we don't claim are completely precise): 

\begin{itemize}
\item Can we axiomatize this notion of ``good" representation ?

\item Can this axiomatization be formulated entirely in the language of symmetric monoidal $C^{*}$-category, or do we need some additional structures on $\Hcal(\Tcal)$ ?

\item Can we characterize the $C^{*}$-categories (possibly with these additional structures) that arise as $\Hcal(\Gcal)$ for some localic stacks $\Gcal$ ?

\item Can we deduce from this a duality between a sub-class of nice ``analytic" localic stacks and certain symmetric monoidal $C^{*}$-categories (possibly with additional structures ?)

\end{itemize}
}

\blockn{As mentioned, the usual Gelfand duality, the duality between commutative monotone complete $C^{*}$-algebras and boolean locales, and the Doplicher-Roberts reconstruction theorem of \cite{doplicher1989new}, can be formulated as special cases of a duality of this form. }

\blockn{We conclude by a few suggestions on the kind of ``additional strucutres" or properties of the $C^{*}$-categories of the form $\Hcal(\Gcal)$ that seems to play an important role.

\begin{itemize}
\item If $a_i$ is a net of operators between two objects $A$ and $B$ of $\Hcal(\Gcal)$ then one can define what it mean for $a_i$ to converge weakly,strongly or strongly-$*$ locally (or ``internally"). In the monotone complete case, normal representations are exactly those who preserve these notions of ``weak" convergence.

\item If $C$ is a small (possibly non-unital) $C^{*}$-category, $F :C \rightarrow \Hcal(\Gcal)$ is any $*$-functor then there is a canonical extension of $F$ from the $C^{*}$-category of $C$-Hilbert modules to $\Hcal(\Gcal)$, defined locally (``internally") as the tensorization of a Hilbert module by the representation ``$F$". This construction seems to be a $C^{*}$-categorical analogue of the classical notions of left/right Kan extension and of weighted (co)-limits, but it does not seems that there is a nice universal property characterizing this notion purely in terms of $C^{*}$-category theory\footnote{But one can give such a universal characterization if one allows the use of the weak convergence mentioned in the first point}. Moreover, the representations classified by $\spec C$ should preserve these weighted (co)limits.
Somehow categories like $\Hcal(\Gcal)$ seems to deserve the name of ``complete" $C^{*}$-category, and normal functor can also be defined as the functors preserving those weighted limits.

\item Every object of $\Hcal(\Gcal)$ has a dual object. Unfortunately, when the object is not locally of finite dimension, there is no clear characterization of what is this dual object in terms of monoidal $C^{*}$-category. Those ``dual" objects seems to play an extremely important role: they are completely central in the proof of the Doplicher-Roberts  reconstruction theorem, and the net of operators $V_u$ (and its dual $V^{*}_U$) used in section \ref{GeomorphonObjec} seem to be a sort of asymptotical replacement for the evaluation and co-evaluation map identifying $l^{2}(X)$ with its own dual and this was a key ingredient in the proof of our theorem.

\end{itemize}

}

\blockn{Hence understanding how we can define what is a ``complete $C^{*}$-category" in the general case and what is the relation of this notion with the notion of ``weak convergence" (these questions are already completely understood in the case of $W^{*}$-categories, see \cite{wstarcat}), and finally understanding what can be says about ``dual objects" in this context seems to be questions that needs to be answered to be able to extend the technique of the present article to general (locally decidable) toposes, and maybe after that to localic groupoids.

Somehow, by working with boolan locally separated toposes in the present paper we did exactly what was needed to avoid those questions: Booleanness ensure that we will only work with monotone complete $C^{*}$-categories, and hence that the question of completeness and normal functors are well understood, and the locale separation hypothesis allows to obtain this net of operators $(V_U)$ which ``rigidify" the notion of dual object at least in $\Hcal^{red}$.}

\bibliography{Biblio}{}
\bibliographystyle{plain}

\end{document}